\documentclass[12pt,DIV=14]{scrartcl}
\usepackage{amsthm,amssymb,amsmath}
\usepackage{enumerate}
\usepackage
[colorlinks=true,citecolor=red,urlcolor=purple,linkcolor=blue,pdfborder={0 0 0}]{hyperref} % enables hyperlinks
\usepackage{tikz}
\usetikzlibrary{calc}

\allowdisplaybreaks[2]

\newtheorem{theorem}{Theorem}[section]
\newtheorem{corollary}[theorem]{Corollary}
\newtheorem{lemma}[theorem]{Lemma}

\theoremstyle{definition}

\theoremstyle{remark}
\newtheorem{remark}[theorem]{Remark}

\providecommand{\C}{\mathbb{C}}

\providecommand{\N}{\mathbb{N}}

\providecommand{\R}{\mathbb{R}}
\providecommand{\Z}{\mathbb{Z}}

\providecommand{\supp}{\operatorname{supp}}
\let\Re\relax\DeclareMathOperator{\Re}{Re}

\providecommand{\Om}{\Omega}
\providecommand{\pOm}{{\partial\Omega}}

\providecommand{\vecr}{\hat\vecx}
\providecommand{\vecx}{\boldsymbol{x}}

\providecommand{\vecphi}{\boldsymbol{\varphi}}
\providecommand{\vecnu}{\boldsymbol{\nu}}

\providecommand{\T}{\mathcal{T}}
\providecommand{\Th}{\mathcal{T}_h}

\providecommand{\TintBR}{{\tilde{T}}}
\providecommand{\TextBR}{{\breve{T}}}
\providecommand{\ThintOm}{{\Th^-}}
\providecommand{\Thstrip}{{\Th^+}}
\providecommand{\ThintBR}{{\widetilde{\T}_h}}
\providecommand{\ThextBR}{{\breve{\T}_h}}

\providecommand{\fcal}{\mathcal{F}}
\providecommand{\ncal}{\mathcal{N}}
\providecommand{\pcal}{\mathcal{P}}

\renewenvironment{thebibliography}[1]{
\begin{oldthebibliography}{#1}
\setlength{\itemsep}{0em}
\setlength{\parskip}{0em}
}
{
\end{oldthebibliography}
}
\title{Finite element solution of a radiation/propagation problem for a Helmholtz equation
with a compactly supported nonlinearity}
\author{Lutz Angermann\thanks{%
Dept.~of Mathematics, Clausthal University of Technology, Erzstr.~1, D-38678 Clausthal-Zellerfeld, Germany,
lutz.angermann@tu-clausthal.de}}
\date{\today}
\begin{document}

\maketitle

\begin{abstract}
A finite element approach for approximating the solution of a mathematical model
for the response of a penetrable, bounded object (obstacle) to the excitation by 
an external electromagnetic field is presented and investigated.
The model consists of a nonlinear Helmholtz equation that is reduced
to a spherical domain.

As a specific example, we consider a finite element method consisting of Courant-type elements
with curved edges at the boundary of a circular computational domain in the two-dimensional case.
We examine this method and more general conforming methods -- including three-dimensional ones --
with comparable properties for their well-posedness; in particular, the validity of a discrete
inf-sup condition of the modified sesquilinear form uniformly with respect to 
both the truncation and the mesh parameters is shown.
Under suitable assumptions to the nonlinearities, a quasi-optimal error estimate is obtained.
Finally, the satisfiability of the approximation property of the finite element space
required for the solvability of a class of adjoint linear problems is discussed. 
\end{abstract}



\setcounter{section}{-1}

\section{Introduction}
The present work deals with the numerical approximation of the solution
of a mathematical model for the response of a penetrable, bounded object (obstacle)
to the excitation by an external electromagnetic field.
Such models are not only interesting from a mathematical point of view,
but also from a physical or engineering user perspective,
for example in the case of optical switches or frequency multipliers.

The paper can be considered as a continuation of the work \cite{Angermann:23f}.
In the latter, the original full-space transmission problem for a nonlinear Helmholtz equation
is transformed into an equivalent boundary value problem on a bounded domain
using a nonlocal Dirichlet-to-Neumann (DtN) operator and analyzed analytically.
It is shown that the transformed nonlinear problem is equivalent to the original problem
and can be uniquely solved under certain conditions.
In addition, the effect of the truncation of the DtN operator on the resulting solution
is examined.
It is shown that the corresponding sesquilinear form satisfies
a parameter-uniform inf-sup condition, that the modified nonlinear problem
has a unique solution, and that the solution error caused by the truncated DtN operator
can be estimated in dependence on the truncation parameter.

Here the question is examined how the solution of the transformed nonlinear problem
can be approximated by means of a finite element method (FEM).
The basic approach and the key steps of the analysis are described using
a relatively simple method in the two-dimensional case -- practically a Courant FEM
with curved geometric elements at the boundary of the circular computational domain.
Although it is well-known that this method is not sufficiently robust with regard
to large wave numbers \cite{Babuska:97b} and that other methods,
such as discontinuous Galerkin methods (e.g., \cite{melenk:13}, \cite{Congreve:19}
to mention a few), are therefore more suitable,
this weak point was deliberately accepted here in order not to let
the presentation swell too much in terms of scope.
A transfer to conforming higher-order FEM is easily possible provided
a uniform discrete inf-sup condition of the truncated sesquilinear form and
suitable approximation capabilities of the approximation space
for the solutions of the original as well as the adjoint linear problem are present.
A discussion of such aspects can be found in the literature
(e.g., \cite{Melenk:10} and some later contributions by these authors).
However, if discontinuous methods are to be used, further investigations are required.

This paper is structured as follows. In the first two sections,
the original problem and the resulting weak formulation,
which forms the starting point for the discretization, are introduced.
More detailed explanations about the model can be found in the work \cite{Angermann:23f}
already mentioned.
The third section is dedicated to the construction of an exemplary finite element space.
The mesh conception consists of using triangle-like elements that have
at most one curved edge at the boundary of the circular computational domain,
and interface-adapted triangles inside the computational domain.
Piecewise polynomial, continuous functions with local polynomial degree one
are then defined over these geometric elements.
For this concept, the estimates of the interpolation error required later on
are derived using a simple continuation technique.
At the end of the section, all the key properties of the finite element discretization
are summarized, as the subsequent investigations also apply to more general conforming methods,
including three-dimensional ones.

The main Section 4 essentially contains the well-posedness investigations
for the finite element discretization.
The results are not limited to the exemplary method described, but are also valid for other,
including three-dimensional, finite element methods,
provided they have the properties described at the end of Section 3.
Under the assumption that the finite element approximation space has appropriate approximation capabilities
to the solutions of a class of adjoint linear problems, the unique solvability
of the linear original problem as well as the uniform
(both with respect to the truncation as well as the mesh parameters) validity of
a discrete inf-sup condition of the truncated sesquilinear form are demonstrated.
After that, the unique solvability of the nonlinear problem -- of course under
appropriate assumptions about the nonlinearities -- is proven.
At the end of the section an error estimate for the finite element solution
of the nonlinear problem is given.
This estimate, together with the estimate of the truncation error from the work \cite[Thm.~27]{Angermann:23f},
yields an estimate of the total error, i.e., the error between the FE approximation
and the solution of the full-space transmission problem.
In the last subsection the satisfiability of the mentioned approximation property
of the finite element space for the solution of a class of adjoint linear problems is discussed.

\section{Problem formulation}
\label{sec:problem}
\begin{figure}[htb]
\begin{center}
\begin{tikzpicture}
\shade [right color=red,opacity=.7] plot [smooth cycle] coordinates {(0,0) (1,1) (2.4,1) (2.3,-1)};
\shade[ball color = blue, opacity = 0.15] (1.5,0) circle [radius = 2];
\draw[->] (-.8,2.1) -- (1.2,1.1);
\node at (2,0) {$\Om$};
\node at (3.8,-1) {$S_R$};
\node at (-.7,1.6) {$u^\mathrm{inc}$};
\end{tikzpicture}
\end{center}
\caption{%
The nonlinear medium $\Om$ is excited by an incident field $u^\mathrm{inc}$}
\label{fig:radscat01}
\end{figure}
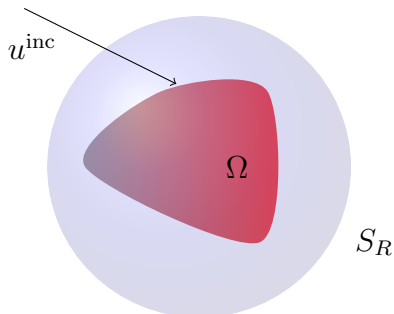

Let $\Om\subset\R^d$ be a bounded domain with a Lipschitz boundary $\pOm$.
It represents a medium that responds nonlinearly to interaction with electromagnetic fields.
Since $\Om$ is bounded, we can choose an open Euclidean $d$-ball $B_R\subset\R^d$ of radius
$R > R_0:=\sup_{\vecx\in\Om}|\vecx| > 0$
with center in the origin such that $\Om\subset B_R$.
The complements of $\Om$ and $B_R$ are denoted by
$\Om^\mathrm{c} := \R^d\setminus\Om$
$B_R^\mathrm{c} := \R^d\setminus B_R$, resp.,
the open complement of $B_R$ is denoted by $B_R^+ := \R^d\setminus\overline{B_R}$
(the overbar over sets denotes their closure in $\R^d$),
and the boundary of $B_R$, the sphere, by $S_R := \partial B_R$ (cf.\ Fig.~\ref{fig:radscat01}).
The open complement of $\Om$ is denoted by $\Om^+ := \R^d\setminus\overline{\Om}$.
By $\vecnu$ we denote the outward-pointing (w.r.t.\ either $\Om$ or $B_R$) unit normal vector
on $\pOm$ or $S_R$, respectively.

Trace operators will be denoted by one and the same symbol $\gamma$; the concrete
meaning (e.g., traces on the common interface of an interior and exterior domain)
will be clear from the context.

The classical direct problem of radiation and propagation of an electromagnetic field
-- actually just one component of it --
by/in the penetrable obstacle $\Om$ is governed by a nonlinear Helmholtz equation
with a variable complex-valued wave coefficient:
\begin{equation}\label{eq:genscalproblem}
-\Delta u(\vecx) - \kappa^2 c(\vecx,u) \,u = f(\vecx,u)
\quad\text{for (almost) all }
\vecx\in\R^d,
\end{equation}
where the wave number $\kappa>0$ is fixed.
The physical properties of the obstacle $\Om$ are described by the coefficient
$c:\;\R^d\times\C\to\C$ (physically the square of the \emph{refractive index}) and
the right-hand side $f:\;\R^d\times\C\to\C$.
In general, both functions are nonlinear and have the following properties:
\begin{equation}\label{eq:ass_nonlin1}
\supp(1-c(\cdot,w)) = \overline\Om
\quad\text{and}\quad
\supp f(\cdot,w) \subset\overline\Om
\quad\text{for all }
w\in\C.
\end{equation}
The function $1-c$ is often called the \emph{contrast function}.
Basically we assume that $c$ and $f$ are Carath\'{e}odory functions,
i.e.,
the mapping $\vecx\mapsto c(\vecx,v)$ is (Lebesgue-)measurable for all $v\in\C$,
and
the mapping $v\mapsto c(\vecx,v)$ is continuous for almost all $\vecx\in\R^d$.
These two conditions imply that $\vecx\mapsto c(\vecx,v(\vecx))$
is measurable for any measurable $v$. The same applies to $f$.

The unknown \emph{total field} $u:\;\R^d\to\C$ should have the following structure:
\begin{equation}\label{eq:solstructure}
u=\begin{cases}
u^\mathrm{rad} + u^\mathrm{inc} & \text{in }\Om^\mathrm{c},\\
u^\mathrm{trans} & \text{in }\Om,
\end{cases}
\end{equation}
where
$u^\mathrm{rad}:\;\Om^\mathrm{c}\to\C$ is the unknown radiated/scattered field,
$u^\mathrm{trans}:\;\Om\to\C$ denotes the unknown transmitted field,
and the incident field $u^\mathrm{inc}\in H^1_\mathrm{loc}(\Om^+)$ is given,
where
\[
H^1_\mathrm{loc}(\Om^+):=\left\{v\in L_{2,\mathrm{loc}}(\Om^+):
\;\varphi v \in H^1(\Om^+)
\text{ for all }
\varphi\in C^\infty_0(\Om^\mathrm{c})
\right\}.
\]
The incident field is usually a (weak) full-space solution of either the homogeneous or
inhomogeneous Helmholtz equation.
Typically it is generated either by concentrated sources located
in a bounded subdomain of $\Om^+$ or by sources at infinity, e.g.\ travelling waves.

The radiated/scattered field $u^\mathrm{rad}$
should satisfy an additional condition, the so-called \emph{Sommerfeld radiation condition}:
\begin{equation}\label{eq:sommerfeld1}
\lim_{|\vecx|\to\infty} |\vecx|^{(d-1)/2}
\left(\vecr\cdot\nabla u^\mathrm{rad} - i\kappa u^\mathrm{rad}\right) = 0
\end{equation}
uniformly for all directions $\vecr:=\vecx/|\vecx|$, where
$\vecr\cdot\nabla u^\mathrm{rad}$ denotes the derivative of $u^\mathrm{rad}$
in radial direction $\vecr$,
cf.\ \cite[eq.~(3.7) for $d=3$, eq.~(3.96) for $d=2$]{Colton:13b}.

\section{Weak formulations}
A weak formulation of the problem \eqref{eq:genscalproblem}--\eqref{eq:sommerfeld1}
can be given as follows (for details see \cite{Angermann:23f}).
Introducing the (complex) linear function spaces
\begin{align*}
H^1_\mathrm{comp}(\Om^+) &:= \left\{v\in H^1(\Om^+):
\;\supp v\ \text{ is compact}\right\},\\
V_{\R^d} &:= \{v\in L_2(\R^d):\; v|_\Om\in H^1(\Om),\ v|_{\Om^+}\in H^1_\mathrm{loc}(\Om^+):
\;\gamma v|_\Om = \gamma v|_{\Om^+} \text{ on }\pOm\},\\
V_{\R^d}^\circ &:= \{v\in L_2(\R^d):\; v|_\Om\in H^1(\Om),\ v|_{\Om^+}\in H^1_\mathrm{comp}(\Om^+):
\;\gamma v|_\Om = \gamma v|_{\Om^+} \text{ on }\pOm\}
\end{align*}
(note the comment at the beginning of Section \ref{sec:problem} on the notation for trace operators),
we say that $u\in V_{R^d}$ is a \emph{weak solution}
to the problem \eqref{eq:genscalproblem}--\eqref{eq:sommerfeld1}
for given $u^\mathrm{inc}\in H^1_\mathrm{loc}(\Om^+)$ if
it has the structure \eqref{eq:solstructure} and
the variational equation
\begin{equation}\label{eq:weakfullspace}
(\nabla u,\nabla v)_\Om + (\nabla u,\nabla v)_{\Om^+} - \kappa^2(c(\cdot,u) u,v)_{\R^d}
=(f(\cdot,u),v)_{\R^d}
\quad\text{for all }
v\in V_{\R^d}^\circ
\end{equation}
together with the Sommerfeld radiation condition \eqref{eq:sommerfeld1}
are satisfied.
Here and in what follows we use a customary notation for various inner products,
that means, for any domain $M\subset\R^d$ with boundary $\partial M$
and appropriately defined functions on $M$ or $\partial M$, we write
\begin{gather*}
(\nabla w,\nabla v)_M := \int_M \nabla w\cdot\nabla\overline{v} d\vecx,\\
(w,v)_M := \int_M w\overline{v} d\vecx,
\qquad
(w,v)_{\partial M} := \int_{\partial M} w\overline{v} ds(\vecx)
\end{gather*}
(the overbar over functions denotes complex conjugation).

Using the so-called DtN technique, the global (full-space) problem
\eqref{eq:weakfullspace} can be transformed into a problem on a bounded
computational domain, here the ball $B_R$.
In the two-dimensional case, the DtN operator $T_\kappa:\; H^{s+1/2}(S_R)\to H^{s-1/2}(S_R)$,
$s\ge 0$, is defined as follows:
\begin{equation}\label{def:2dDtNb}
\begin{aligned}
T_\kappa u(\vecx)
&:= \frac{1}{R}\sum_{n\in\Z} Z_n(\kappa R)u_n(R)Y_n(\hat\vecx),
\quad
\vecx=R\hat\vecx\in S_R,
\end{aligned}
\end{equation}
where
$\ Z_n(\xi):= \xi H_n^{(1)'}(\xi)/H_n^{(1)}(\xi)$,
$H_n^{(1)}$ are the cylindrical Hankel functions of the first kind of order $n$
\cite[Sect.~10.2]{NIST:23},
$Y_n$ are the circular harmonics defined by
$\ Y_n(\varphi) := (2\pi)^{-1/2}e^{in\varphi}, \ n\in\Z$,
(identifying $Y_n(\hat\vecx)$ with $Y_n(\varphi)$ for
$\vecx=r\hat\vecx=r(\cos\varphi,\sin\varphi)^\top$
in polar coordinates),
and $u_n(R)$ are the Fourier coefficients of $u|_{S_R}$:
\begin{equation}\label{eq:Fouriercoeff2d}
u_n(R) := (u(R\cdot),Y_n)_{S_1} = \int_{S_1}u(R\hat\vecx)\overline{Y_n}(\hat\vecx)ds(\hat\vecx)
=\int_0^{2\pi} u(R,\varphi)\overline{Y_n}(\varphi)d\varphi
\end{equation}
(here we identify $u(\vecx)$ with $u(r,\varphi)$;
$ds(\hat\vecx)$ denotes the Lebesgue arc length element).
For $d=3$ we define
\begin{equation}\label{def:3dDtNb}
\begin{aligned}
T_\kappa u(\vecx)
&:= \frac{1}{R}\sum_{n\in\N_0}\sum_{|m|\le n}
z_n(\kappa R)u_n^m(R)Y_n^m(\hat\vecx),
\quad
\vecx=R\hat\vecx\in S_R,
\end{aligned}
\end{equation}
where
$z_n(\xi):= \xi h_n^{(1)'}(\xi)/h_n^{(1)}(\xi)$,
$h_n^{(1)}$ are the spherical Hankel functions of the first kind of order $n$
\cite[Sect.~10.47]{NIST:23},
$Y_n^m$ are the spherical harmonics defined by
\[
Y_n^m(\varphi,\theta) := \sqrt{\frac{2n+1}{4\pi}\,\frac{(n-|m|)!}{(n+|m|)!}}\,
P_n^{|m|}(\cos\theta)e^{im\varphi},
\quad
n\in\N_0,\ |m|\le n,
\]
(identifying $Y_n^m(\hat\vecx)$ with $Y_n^m(\varphi,\theta)$ for
$\hat\vecx=(\cos\varphi\sin\theta,\sin\varphi\sin\theta,\cos\theta)^\top$
in spherical coordinates),
where
$P_n^m$ are the associated Legendre functions of the first kind \cite[Sect.~14.21]{NIST:23},
and $u_n^m(R)$ are the Fourier coefficients of $u|_{S_R}$:
\begin{equation}\label{eq:Fouriercoeff3d}
\begin{aligned}
u_n^m(R) = (u(R\cdot),Y_n^m)_{S_1}
&= \int_{S_1}u(R\hat\vecx)\overline{Y_n^m}(\hat\vecx)ds(\hat\vecx)\\
&=\int_0^{2\pi}\int_0^\pi u(R,\varphi,\theta)
\overline{Y_n^m}(\varphi,\theta)\sin\theta d\theta d\varphi
\end{aligned}
\end{equation}
(here we identify $u(\vecx)$ with $u(r,\varphi,\theta)$;
$ds(\hat\vecx)$ is the Lebesgue surface area element).

In \cite[Lemma~7]{Angermann:23f} it was shown that the weak formulation
\eqref{eq:weakfullspace} is equivalent to the following problem on the ball $B_R$:

Find $u\in V$ such that
\begin{equation}\label{eq:weakball}
\begin{aligned}
& (\nabla u,\nabla v)_\Om + (\nabla u,\nabla v)_{B_R\setminus\overline\Om}
- \kappa^2(c(\cdot,u) u,v)_{B_R} - (T_\kappa u,v)_{S_R}\\
&= (f(\cdot,u),v)_{B_R}
- (T_\kappa u^\mathrm{inc},v)_{S_R} + (\vecr\cdot\nabla u^\mathrm{inc},v)_{S_R}
\quad\text{for all }
v\in V,
\end{aligned}
\end{equation}
where
\[
V := \{v\in L_2(B_R):\; v|_\Om\in H^1(\Om),\ v|_{B_R\setminus\overline\Om}\in H^1(B_R\setminus\overline\Om):
\;\gamma v|_\Om = \gamma v|_{B_R\setminus\overline\Om} \text{ on }\pOm\}.
\]
\begin{remark}\label{rem:V=H1}
In fact, $V = H^1(B_R)$. Hence the first two terms on the left-hand side in \eqref{eq:weakball}
can be replaced by $(\nabla u,\nabla v)_{B_R}$.
\end{remark}
\begin{proof}
We demonstrate the $V \subset H^1(B_R)$-part only.
So let $v\in V$ be fixed and $\vecphi \in C_0^\infty(B_R)^d$ be arbitrary.
Then
\begin{align*}
\int_{B_R} v\,\nabla\cdot\vecphi \, dx
&= \int_\Om v\,\nabla\cdot\vecphi \, dx + \int_{B_R\setminus\overline\Om} v\,\nabla\cdot\vecphi \, dx\\
&= - \int_\Om \nabla(v|_\Om)\cdot\vecphi \, dx + \int_\pOm \gamma v|_\Om \, \vecnu\cdot\vecphi \, ds(\vecx)\\
& \qquad
- \int_{B_R\setminus\overline\Om} \nabla(v|_{B_R\setminus\overline\Om})\cdot\vecphi \, dx
+ \int_{\partial(B_R\setminus\overline\Om)} \gamma v|_{B_R\setminus\overline\Om} \, \vecnu\cdot\vecphi \, ds(\vecx).
\end{align*}
Defining $\nabla v:\; B_R \to \C$ by
\[
\nabla v := \begin{cases}
\nabla(v|_\Om) &\text{on }\Om,\\
\nabla(v|_{B_R\setminus\overline\Om}) &\text{on }B_R\setminus\overline\Om,
\end{cases}
\]
we get
\begin{align*}
\int_{B_R} v\,\nabla\cdot\vecphi \, dx
&= - \int_{B_R} \nabla v \cdot\vecphi \, dx + \int_{S_R} v \, \vecnu\cdot\vecphi \, ds(\vecx)
+ \int_\pOm [\gamma v|_\Om - \gamma v|_{B_R\setminus\overline\Om}] \, \vecnu\cdot\vecphi \, ds(\vecx).
\end{align*}
The second integral on the right-hand side vanishes due to the zero boundary values of $\vecphi$
on $S_R$, whereas the third integral vanishes thanks to the zero jump of $\gamma v$ on $\pOm$.
Thus $\nabla v$ is the weak gradient of $v$ on $B_R$ in the sense of distributions.
It is an element of $L_2(B_R)^d$ since its restrictions to $\Om$ and $B_R\setminus\overline\Om$
belong to $L_2(B_R)^d$. Hence $v\in H^1(B_R)$.
\end{proof}

On the space $V$, we use the standard seminorm and norm:
\begin{align*}
|v|_V &:= \left(\|\nabla v\|_{0,2,\Om}^2 + \|\nabla v\|_{0,2,B_R\setminus\overline\Om}^2\right)^{1/2}
= \|\nabla v\|_{0,2,B_R},\\
\|v\|_V &:= \left(|v|_V^2 + \|v\|_{0,2,B_R}^2\right)^{1/2} = \|v\|_{1,2,B_R}.
\end{align*}
For $\kappa>0$, the following so-called \emph{wavenumber dependent norm} on $V$
is also common:
\begin{equation}\label{eq:wavenumbernorm}
\|v\|_{V,\kappa} := \left(|v|_V^2 + \kappa^2\|v\|_{0,2,B_R}^2\right)^{1/2}.
\end{equation}
Obviously, the standard norm and the wavenumber dependent norm
are equivalent on $V$, i.e., it holds
\begin{equation}\label{eq:wavenumbernormequiv}
C_- \|v\|_V \le \|v\|_{V,\kappa} \le C_+ \|v\|_V
\quad\text{for all } v\in V
\end{equation}
with
$C_- := \min\{1;\kappa\}$ and $C_+ := \max\{1;\kappa\}$.

Unfortunately, this problem is still difficult to solve numerically
since the DtN operators $T_\kappa$ are nonlocal.
This leads to the idea of truncating the representations
\eqref{def:2dDtNb}, \eqref{def:3dDtNb}
as follows. For $N\in\N_0$, the \emph{truncated DtN operators} are defined by
\begin{align}
T_{\kappa,N} u(\vecx)
&:= \frac{1}{R}\sum_{|n|\le N} Z_n(\kappa R)u_n(R)Y_n(\hat\vecx),
\quad
\vecx=R\hat\vecx\in S_R\subset\R^2,
\label{def:2dDtNtrunc}\\
T_{\kappa,N} u(\vecx)
&:= \frac{1}{R}\sum_{n=0}^N\sum_{|m|\le n}
z_n(\kappa R)u_n^m(R)Y_n^m(\hat\vecx),
\quad
\vecx=R\hat\vecx\in S_R\subset\R^3.
\label{def:3dDtNtrunc}
\end{align}

If we define, for all $w,v\in V$,
\begin{align}
a_N(w,v) &:= (\nabla w,\nabla v)_{B_R}
- \kappa^2 (w,v)_{B_R} - (T_{\kappa,N} w,v)_{S_R},\nonumber\\
n_N(w,v) &:= \langle\ell^\mathrm{contr}(w),v\rangle + \langle\ell^\mathrm{src}(w),v\rangle
+\langle\ell_N^\mathrm{inc},v\rangle
\label{def:n_N}\\
&:= \kappa^2((c(\cdot,w) - 1) w,v)_{B_R} + (f(\cdot,w),v)_{B_R}
- (T_{\kappa,N} u^\mathrm{inc},v)_{S_R} + (\vecr\cdot\nabla u^\mathrm{inc},v)_{S_R},
\nonumber
\end{align}
we can formulate the following problem:

Find $u_N\in V$ such that
\begin{equation}\label{eq:intscalproblem4}
a_N(u_N,v) = n_N(u_N,v)
\quad\text{for all }
v\in V.
\end{equation}
Under the assumption that $R \ge R_0$ and $\kappa \ge \kappa_0 >0$ with $\kappa_0 \ge 1$ if $d=2$,
it can be shown that the truncated sesquilinear form $a_N$
satisfies a uniform inf-sup condition in the following sense \cite[Lemma~25]{Angermann:23f}:

There exists a number $N^\ast\in\N$ such that
\begin{equation}\label{eq:infsuptrunc}
\beta_{N^\ast}(R,\kappa) := \inf_{w\in V\setminus\{0\}}\sup_{v\in V\setminus\{0\}}
\frac{|a_N(w,v)|}{\|w\|_{V,\kappa} \|v\|_{V,\kappa}}>0
\end{equation}
is independent of $N\ge N^\ast$.
Moreover, the inf-sup constant $\beta_{N^\ast}(R,\kappa)$ can be expressed
in terms of the inf-sup constant $\beta(R,\kappa)$ of the sesquilinear form $a$.

This property guarantees the stable solvability of the problem \eqref{eq:intscalproblem4}
in the case that the right-hand side is a linear continuous functional over V,
and enables to prove the following theorem \cite[Thm.~26]{Angermann:23f}.
\begin{theorem}\label{th:truncscalee}
Under the above assumptions to $R$ and $\kappa$,
let the functions $c$ and $f$ generate locally Lipschitz continuous
Nemycki operators in $V$ and assume that there exist functions $w_f,w_c\in V$ such that
$f(\cdot,w_f)\in L_{p_f/(p_f-1)}(\Om)$ and $c(\cdot,w_c)\in L_{p_c/(p_c-2)}(\Om)$,
respectively.

Furthermore let
$u^\mathrm{inc}\in H^1_\mathrm{loc}(\Om^+)$ be such that
additionally $\Delta u^\mathrm{inc}\in L_{2,\mathrm{loc}}(\Om^+)$ holds.

If there exist numbers $\varrho>0$ and $L_\mathcal{F}\in (0,\beta_{N^\ast}(R,\kappa))$
such that the following two conditions
\begin{align}
&\kappa^2\left[\|c(\cdot,w_c)-1\|_{0,\tilde{q}_c,\Om}
+ \|L_c(\cdot,w,w_c)\|_{0,q_c,\Om}(\varrho + \|w_c\|_V)\right]\varrho\nonumber\\
&\quad + \left[\|f(\cdot,w_f)\|_{0,\tilde{q}_f,\Om}
+ \|L_f(\cdot,w,w_f)\|_{0,q_f,\Om}(\varrho + \|w_f\|_V)\right]\label{eq:thtrunceu1}\\
&\quad
+ C_\mathrm{tr}\|\vecr\cdot\nabla u^\mathrm{inc} - T_{\kappa,N} u^\mathrm{inc}\|_{-1/2,2,S_R}
\le \varrho\beta_{N^\ast}(R,\kappa),\nonumber\\
&\kappa^2\left[\|L_c(\cdot,w,v)\|_{0,q_c,\Om}\varrho
+ \|c(\cdot,w_c)-1\|_{0,\tilde{q}_c,\Om}
+ \|L_c(\cdot,w,w_c)\|_{0,q_c,\Om}(\varrho + \|w_c\|_V)\right]\nonumber\\
&\quad + \|L_f(\cdot,w,v)\|_{0,q_f,\Om} \le L_\mathcal{F}
\label{eq:thtrunceu2}
\end{align}
are satisfied for all $w,v\in K_\varrho^\mathrm{cl}:=\{v \in V:\; \|v\|_V \le \varrho\}$,
then the problem \eqref{eq:intscalproblem4} has a unique solution
$u_N \in K_\varrho^\mathrm{cl}$ for all $N \ge N^\ast$.
\end{theorem}
Finally, we have proved the following result \cite[Thm.~27]{Angermann:23f}.
\begin{theorem}\label{th:truncscalerrest}
Let the assumptions of the above Thm.~\ref{th:truncscalee}
with respect to $R$, $\kappa$, $c$, and $f$ be satisfied.
Then, if the Lipschitz constant $L_\mathcal{F}$ is sufficiently small,
i.e., satisfies
\[
L_\mathcal{F}\le\min\left\{\beta(R,\kappa),\beta_{N^\ast}(R,\kappa),
\frac{C_-^2}{4\tilde{\varrho}\big(1+2\kappa^2 C_-^{-1}C(R,\kappa) C_\mathrm{emb}\big)}
\right\},
\]
there exists a constant $c>0$ independent of $N\ge N^\ast$
such that the following estimate holds:
\[
c\|u-u_N\|_V
\le c_+(N,u) \|u\|_V
+ c_+(N,u^\mathrm{inc}) \|u^\mathrm{inc}\|_V\\
+ \frac{1}{(1+N^2)^{1/2}} \|u_N\|_V.
\]
\end{theorem}

\section{Discretization}
To approximate the solution of \eqref{eq:intscalproblem4}, conforming finite element methods will be used.
We describe the details of an exemplary method below for the two-dimensional case.
For the three-dimensional case, we refer to the literature, e.g., \cite{Bernardi:89}.
There are basically two additional problems to be dealt with compared
to classical methods on polyhedral domains.
The first concerns the fact that the boundary $S_R$ of the computational domain
is a sphere, i.e., a smooth manifold that cannot be represented exactly by polynomial patches.
The latter would be a prerequisite for using conventional isoparametric finite elements.
An alternative way could be an isogeometric approach, where the sphere can be modeled exactly
by a NURBS surface \cite{Cobb:88}.
However, this approach does not have to be applied in our procedure
since we can use an explicit parameterization of the sphere.
The second issue concerns the transmission problem, more precisely
the treatment of the interface $\pOm$.
Here we apply a so-called \emph{interface-fitted method}, in which the mesh is adapted to the interface.
A different class of methods is formed by interface-unfitted methods
whose meshes ignore the interface to some extend.
A typical representative is the so-called CutFEM, see, e.g., \cite{Burman:15}.

Formally, the problem of designing finite element methods on domains with curved boundaries
is fairly well-understood nowadays, see, e.g., \cite{Kawecki:20}.
However, since fundamental aspects of discretization are to be discussed here first,
we restrict ourselves to a simpler, less general approach with regard to geometric partitioning.

\subsection{Geometric discretization ($d=2$)}

The admissible partitions of the open disk $B_R$ are described via a three-step procedure
to obtain a so-called \emph{interface-fitted partition} $\Th$ of $B_R$.
\begin{enumerate}[(T1)]
\item\label{item:T1}
We start with a polygonal approximation $\Om_h$ of $\Om$,
where $\Om_h$ is to be understood as an open set, too.
To do this, let $\ThintOm$ be a finite collection of pairwise disjoint, open triangles,
generally denoted by $T$, such that
\begin{enumerate}[(i)]
\item
$\forall T\in\ThintOm:\;T\cap\Om\ne\emptyset$ (\emph{position property}).
\item
$\Om_h$ is the interior of
\[
\bigcup_{T\in\ThintOm}\overline{T}.
\]
\item
$\forall T,T'\in\ThintOm$: The set $\overline{T}\cap\overline{T}'$
is either empty, a common vertex or a common edge of $T$ and $T'$
(\emph{matching property, consistency} or \emph{conformity}).
\item
All the vertices of the polygon $\Om_h$ lie on $\pOm$ (\emph{fitting property}).
\end{enumerate}
The mesh parameter of this triangulation is defined as
$h := \max_{T\in\ThintOm} h_T$, where $h_T$ is the diameter of $T$.
\item\label{item:T2}
In the second step we design a triangulation $\Thstrip$
of a polygonal subdomain of $B_R\setminus\overline{\Om}_h$ and combine it
with $\ThintOm$ to form a triangulation of a polygonal subdomain of $B_R$ as follows.

We take an approximating polygon $B_{Rh}$ of the disk $B_R$
with mesh parameter $h$ of the same magnitude as that of $\ThintOm$.
This can be obtained by means of the same approach as described for $\Om$,
but replacing $\Om$ by the disk $B_R$
(however, we do not need this triangulation in its entirety for the following).
We also require that each triangle has at most one edge on $\partial B_{Rh}$.
This is not an essential restriction since we will require later that the
final triangulation is shape-regular, which (at least for sufficiently large $R$)
excludes too flat triangles.
Since $B_R$ is a convex set, $B_{Rh}\subset B_R$.
Then we triangulate the set $B_{Rh}\setminus\overline\Om_h$
by a triangulation $\Thstrip$ in such a way that the combined triangulation
$\ThintBR:=\ThintOm\cup\Thstrip$ of $B_{Rh}$ possesses the matching property
analogous to (T\ref{item:T1})(iii).

\item\label{item:T3}
The final step is to replace all the boundary triangles with triangles
that have a curved edge.

Denote by $\ThintBR^\fcal$ the subset of $\ThintBR$
consisting of triangles having one edge on $\partial B_{Rh}$.

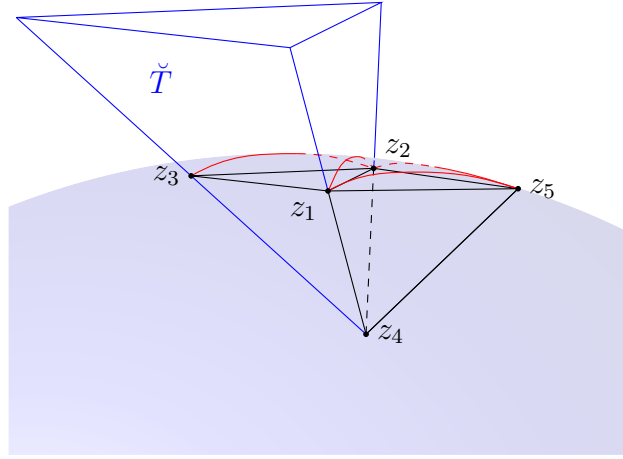
\begin{figure}[h]
\begin{center}
\begin{tikzpicture}
\path[clip](-3.5,6) rectangle(4.9,12.6);
\coordinate (O) at (0,0);
\shade[ball color = blue, opacity = 0.15] (O) circle [radius = 10];

\coordinate (z1) at (3,9.53);
\coordinate (z2) at (1.1,9.94);
\coordinate (z3) at (1.5,7.6);
\coordinate (z4) at (-1.3,9.9);
\coordinate (z5) at (0.2,8.2);

\coordinate (z1b) at ($2*(z1)-(z3)$);
\coordinate (z2b) at ($2*(z2)-(z3)$);
\coordinate (z2c) at ($2*(z2)-(z5)$);
\coordinate (z4b) at ($2*(z4)-(z5)$);

\coordinate (z12) at (0.89,9.955);
\coordinate (z42) at (0.1,9.995);
\coordinate (z25) at (1.9,9.82);
\coordinate (z15) at (1.75,9.75);

\draw[blue] (z1) -- (z1b) -- (z2b) -- (z2);
\draw[blue] (z2) -- (z2c) -- (z4b) -- (z4);

\draw[fill=red!20] (z3) -- ($(O) + (72.5:10)$) arc (72.5:83.7:10) -- cycle;
\draw[fill=red!20] (z5) -- ($(O) + (83.7:10)$) arc (83.7:97.5:10) -- cycle;

\draw (z1) -- (z2) -- (z3) -- cycle;
\draw (z2) -- (z4) -- (z5) -- cycle;

\draw[fill,black] (z1) circle (0.03);
\draw[fill,black] (z2) circle (0.03);
\draw[fill,black] (z3) circle (0.03);
\draw[fill,black] (z4) circle (0.03);
\draw[fill,black] (z5) circle (0.03);
\node[below] at (z5) {$x_T^\fcal$};

\node[blue] at (-1.7,11) {$\breve{T}$};
\node[red] at (-0.1,9.5) {$T$};
\node at (0.1,9) {$\tilde{T}$};
\node at (-2.5,7) {$B_R$};
\end{tikzpicture}
\end{center}
\caption{Part of the mesh near the boundary $S_R$ with two superelements (blue)}
\label{fig:meshpart}
\end{figure}

We now consider a single triangle $\TintBR\in\ThintBR^\fcal$
with an edge $\tilde{F}\subset\partial B_{Rh}$,
and replace this edge by that curved segment of $S_R$
which has the same vertices as $\tilde{F}$.
In this way, a triangle-like element $T$ with a curved edge results.
By $\Th^\fcal$ we denote the set of all these triangle-like elements with a curved edge,
and $\Th$ is now the partition of $B_R$ resulting from a modification of $\ThintBR$
in which all elements $\TintBR$ belonging to $\ThintBR^\fcal$
are replaced by the corresponding elements of $\Th^\fcal$.
\end{enumerate}

At this point, for later use, we also define geometric elements which we will call
\emph{triangular superelements} or \emph{supertriangles}.
Since each $T\in\Th^\fcal$ has a vertex in the open disk $B_R$, denoted by $x_T^\fcal$,
we can associate with $T$ a supertriangle $\TextBR$, which is obtained by doubling the length
of each edge of $T$ that passes through the vertex $x_T^\fcal$, see Fig.~\ref{fig:meshpart}.
If we replace the elements $T\in\Th^\fcal$ in $\Th$ by their corresponding superelements,
we obtain an overlay triangulation $\ThextBR$.
Accordingly we denote by $\breve{\Om}_h$ the interior of
\[
\bigcup_{T\in\Th}\overline{\TextBR}.
\]
We now assume that the mesh size $h$ of $\Th$
(which is, by definition, equal to the mesh size of $\ThintBR$)
and the radius $R$ are such that
\[
T\subset\TextBR
\quad\text{for all }
T\in\Th^\fcal.
\]
Furthermore we assume that the family of triangulations $\ThintBR$
is shape-regular (i.e., regular in the sense of \cite[Ch.~3.1]{Ciarlet:02b}).

Next we define, for a partition $\Th$ of $B_R$ as defined above,
a discrete function space as follows:
\[
V_h:=\left\{v_h\in H^1(B_R):\, v_h|_T\in\pcal_1(T) \text{ for all } T\in\Th\right\},
\]
where $\pcal_1(T)$ denotes the set of polynomials of maximum degree one on $T$.
Since the elements of $V_h$ are continuous functions on $B_R$, the inclusion
$V_h\subset V$ holds.

As the set of degrees of freedom, i.e., as the set of functionals that uniquely define
a concrete element from $V_h$, we choose node-oriented Lagrangian functionals,
i.e., function evaluations in the nodes of $\Th$.

\subsection{Interpolation}
\label{subsec:interpolation}
After these preparations, we can also introduce a canonical interpolation operator as follows.
If $\overline{\ncal}_{\Th}$ denotes the collection of all vertices of the elements of the partition
$\T_h$, we set
\[
I_h:\;C(\overline{B}_R)\to V_h:\;v\mapsto I_hv:=\sum_{y\in\overline{\ncal}_{\Th}} v(y)\phi_y,
\]
where $\phi_y\in V_h$ is the canonical basis element associated with the node $y\in\overline{\ncal}_{\Th}$,
i.e.
\[
\phi_y(x) = \begin{cases}
1, & x=y,\\
0, & x\in\overline{\ncal}_{\Th}\setminus\{y\},
\end{cases}
\quad
x\in\overline{\ncal}_{\Th}.
\]
In the case of local consideration we can introduce analogous concepts.
Let $\ncal_T$ denote the set of all vertices of the element $T\in\T_h$.
Then
\[
I_T:\;C(\overline{T})\to\pcal_1(T):\;v\mapsto I_Tv:=\sum_{y\in\ncal_T} v(y)\phi_y
\]
is the corresponding local interpolation operator.

Since the basis functions $\phi_y$ are polynomials on $T$,
their domain of definition can be canonically extended to the whole of $\R^2$,
so that in particular they are well-defined on the associated supertriangle $\TextBR$
by restricting their $\R^2$-extension to $\TextBR$.
Therefore, if $\breve{\phi}_y$ denote these extended basis functions,
we can also define a local interpolation operator on $\TextBR$:
\[
I_{\TextBR}:\;C(\overline{\TextBR})\to\pcal_1(\TextBR):\;
\breve{v}\mapsto I_{\TextBR}\breve{v}:=\sum_{y\in\ncal_T} \breve{v}(y)\breve{\phi}_y.
\]
That is, the interpolation operator $I_{\TextBR}$ is defined on the same node set
as the interpolation operator $I_T$.

An error estimate can be obtained in the usual way, where the only difference
from the standard case being the position of the interpolation nodes
for the elements $\breve{T}$ that extend $T\in\Th^\fcal$.
Therefore we describe only this case; the appropriate local estimates
for $T\in\Th\setminus\Th^\fcal$ are standard \cite[Ch.~3.1]{Ciarlet:02b}.
If $\hat{T}$ denotes the reference unit (open) triangle, we can introduce the following
affine-linear reference mapping for $\TextBR$:
\[
F_{\TextBR}:\;\hat{T}\to\TextBR:\;
\hat{x}\mapsto F_{\TextBR}\hat{x}:=B_{\TextBR}\hat{x} + x_T^\fcal,
\quad
B_{\TextBR}\in\R^{2,2},\ \mbox{det}\,B_{\TextBR}\ne 0.
\]
The origins of the interpolation nodes are $\hat{x}_1:=\frac12 e_1$, $\hat{x}_2:=\frac12 e_2$,
and $\hat{x}_3:=0$, where $e_1$, $e_2$ are the canonical unit vectors of $\R^2$.

The interpolation operator on $\hat{T}$ is defined as follows:
\[
I_{\hat{T}}:\;C(\overline{\hat{T}})\to\pcal_1(\hat{T}):\;
\hat{v}\mapsto I_{\hat{T}}\hat{v}:=\sum_{\hat{y}\in\hat{\ncal}} \hat{v}(\hat{y})\hat{\phi}_{\hat{y}}
\]
with $\hat{\ncal}:=\{\hat{x}_j\}_{j=1}^3$ and
$\hat{\phi}_{\hat{y}}(\hat{x}) = \begin{cases}
1, & \hat{x}=\hat{y},\\
0, & \hat{x}\in\hat{\ncal}\setminus\{\hat{y}\},
\end{cases}$
\quad
$\hat{x}\in\hat{\ncal}$.

Note that, for $\hat{v}=\breve{v}\circ F_{\TextBR}$,
\[
I_{\hat{T}}\hat{v}\circ F_{\TextBR}^{-1}
=\sum_{\hat{y}\in\hat{\ncal}} \hat{v}(\hat{y})\hat{\phi}_{\hat{y}}\circ F_{\TextBR}^{-1}
=\sum_{y\in\ncal_T} \breve{v}(y)\breve{\phi}_y
=I_{\TextBR}\breve{v},
\]
and, since a first degree polynomial is uniquely defined by its values on $\hat{\ncal}$,
\begin{equation}\label{eq:P1invariance}
I_{\TextBR}\breve{q} = \breve{q}
\quad\text{for all }
\breve{q}\in\pcal_1(\TextBR).
\end{equation}
Now it is easy to see that the interpolation operator $I_{\hat{T}}$ is continuous if
considered as an operator into a normed space of functions defined on $\hat{T}$,
in particular into $H^r(\hat{T})$ for $r\in\{0,1,2\}$:
\[
\|I_{\hat{T}}\hat{v}\|_{r,2,\hat{T}}
\le \sum_{\hat{y}\in\hat{\ncal}} |\hat{v}(\hat{y})|\|\hat{\phi}_{\hat{y}}\|_{r,2,\hat{T}}
\le \|\hat{v}\|_{0,\infty,\hat{T}}
\sum_{\hat{y}\in\hat{\ncal}} \|\hat{\phi}_{\hat{y}}\|_{r,2,\hat{T}}.
\]
Furthermore, Sobolev's embedding theorem \cite[Thm.~4.12]{Adams:03} yields
\[
\|I_{\hat{T}}\hat{v}\|_{r,2,\hat{T}}
\le \hat{C}_\mathrm{emb}\|\hat{v}\|_{2,2,\hat{T}}
\]
with a constant $\hat{C}_\mathrm{emb}>0$.
Then we can write
\[
\|\hat{v} - I_{\hat{T}}\hat{v}\|_{r,2,\hat{T}}
\le \|\hat{v}\|_{r,2,\hat{T}} + \|I_{\hat{T}}\hat{v}\|_{r,2,\hat{T}}
\le \|\hat{v}\|_{r,2,\hat{T}} + \hat{C}_\mathrm{emb}\|\hat{v}\|_{2,2,\hat{T}}
\le (1 + \hat{C}_\mathrm{emb})\|\hat{v}\|_{2,2,\hat{T}}.
\]
Hence the very left-hand side is a continuous sublinear (i.e., subadditive and homogeneous)
functional on $H^2(\hat{T})$.
Moreover, due to \eqref{eq:P1invariance} it vanishes on $\pcal_1(\TextBR)$, so that
we may apply Bramble-Hilbert's lemma \cite[Ch.~3.1]{Ciarlet:02b},
which actually also applies to sublinear functionals:
\[
\|I_{\hat{T}}\hat{v}\|_{r,2,\hat{T}}
\le C_\mathrm{BH}|\hat{v}|_{2,2,\hat{T}}
\quad\text{for all }
\hat{v}\in H^2(\hat{T}).
\]
The back transformation leads to the estimate
\[
\|\breve{v} - I_{\TextBR}\breve{v}\|_{r,2,\TextBR}
\le C_\mathrm{li} h_{\TextBR}^{2-r}|\breve{v}|_{2,2,\TextBR},
\]
where the constant $C_\mathrm{li}>0$ depends on the shape regularity parameter of the family of $\ThextBR$.
Since the elements of $\ThextBR$ results from a similarity transformation
of the elements of $\ThintBR$, the shape parameters of both triangulations are the same.

Now we are prepared to demonstrate the final global interpolation error estimate.
We make use of the fact that any element $v\in H^2(B_R)$ can be extended to
an element $E_{B_R}v\in H^2(\R^2)$ in a bounded manner, i.e., the extension
operator $E_{B_R}:\;H^2(B_R) \to H^2(\R^2)$ is continuous \cite[Thm.~A.4]{McLean:00}.
Then the restrictions of $E_{B_R}v$ to each $\TextBR$ belong to $H^2(\TextBR)$
and can be used as the elements $\breve{v}$ in the above estimates:
\begin{align*}
\|v - I_h v\|_{r,2,B_R}^2
&= \sum_{T\in\Th} \|v - I_h v\|_{r,2,T}^2
\le \sum_{T\in\Th} \|\breve{v} - I_{\TextBR}\breve{v}\|_{r,2,\TextBR}^2\nonumber\\
&\le C_\mathrm{li}^2 \sum_{T\in\Th} h_{\TextBR}^{2(2-r)}|\breve{v}|_{2,2,\TextBR}^2
\le 8 C_\mathrm{li}^2 h^{2(2-r)} |E_{B_R}v|_{2,2,\breve{\Om}_h}^2\\
&\le 8 C_\mathrm{li}^2 h^{2(2-r)} |E_{B_R}v|_{2,2,\R^2}^2
\le 8 C_\mathrm{li}^2 C_\mathrm{ext}^2 h^{2(2-r)} |v|_{2,2,B_R}^2,
\end{align*}
where we have used the relation $h_{\TextBR}\le 2h$ and the fact that
neighbouring supertriangles can intersect, cf.~Fig.~\ref{fig:meshpart}.

Thus we have proved the following interpolation error estimate.
\begin{theorem}\label{th:interrest}
Let a family $(\Th)_{h\in\mathcal{H}}$ of partitions $\Th$ of $B_R$ be given
which satisfy the asssumptions (T\ref{item:T1})--(T\ref{item:T3}),
where the set $\mathcal{H}$ consists of the members of a sequence of positive real numbers
(mesh parameters) with accumulation point 0.
Then the restriction of the interpolation operator $I_h$ to $H^2(B_R)$,
also denoted by $I_h$, is well-defined and satisfies the following error estimate
for $r\in\{0,1,2\}$:
\[
\|v - I_h v\|_{r,2,B_R} \le C_\mathrm{ip} h^{2-r} |v|_{2,2,B_R}
\quad\text{for all }
v\in H^2(B_R),
\]
where $C_\mathrm{ip}>0$ is a constant independent of $h$, $r$ and $v$.
\end{theorem}

Unfortunately, it is not possible to construct finite element spaces in a similarly simple way
for the case $d=3$ \cite[p.~1356]{Schatz:98b}, and therefore we only refer to other sources here,
e.g., \cite[Sect.~6.3]{Bernardi:89}.
As the discussion in Section \ref{sec:fem} below will show,
only a few essential properties of the finite element spaces are relevant to the results presented there.
We therefore summarize these essential properties here in such a way that they also cover the case $d=3$:

There exists a collection of finite-dimensional subspaces $V_h$ of $V$
such that the following properties are satisfied for $d=2$ or $d=3$:
\begin{enumerate}[(D1)]
\item\label{item:D1}
The parameter $h$ is an element of a sequence of positive real numbers with accumulation point zero.
(Typically it is a parameter of the geometric partitioning of the computational domain.)
\item
$V_h \subset V$ for all $h$ from that sequence.
\item
The dimension of $V_h$ increases without bound if $h$ tends to zero.
\item\label{item:D4}
There exist interpolation operators $I_h:\;H^2(B_R)\to V_h$ that
satisfy the following error estimate for $r\in\{0,1,2\}$:
\[
\|v - I_h v\|_{r,2,B_R} \le C_\mathrm{ip} h^{2-r} |v|_{2,2,B_R}
\quad\text{for all }
v\in H^2(B_R),
\]
where the factor $C_\mathrm{ip}>0$ is independent of $h$, $r$ and $v$.
\end{enumerate}

\section{The finite element method}
\label{sec:fem}
The application of the above setting (D\ref{item:D1}) -- (D\ref{item:D4})
to the truncated problem \eqref{eq:intscalproblem4} leads to the following finite element method:

Find $u_h\in V_h$ such that
\begin{equation}\label{eq:intscalfemproblem4}
a_N(u_h,v_h) = n_N(u_h,v_h)
\quad\text{for all }
v_h\in V_h.
\end{equation}
\subsection{Solvability and stability}
We start with the proof of existence and uniqueness of a solution of the corresponding linear
problem.
To prepare this, we introduce the following auxiliary adjoint problem (cf.\ \cite[p.~43]{McLean:00}):

Given $f\in L_2(B_R)$, find $w_N\in V$ such that
\begin{equation}\label{eq:lindualproblem}
\overline{a_N(v,w_N)} = (v,f)_{B_R}
\quad\text{for all } v\in V.
\end{equation}
According to \cite[Thm.~24]{Angermann:23f} the adjoint problem possesses a unique solution $w_N\in V$.
Therefore we may formally write $w_N=(\mathcal{A}_N^*)^{-1}f$,
where $(\mathcal{A}_N^*)^{-1}:\;L_2(B_R)\to V$ denotes the solution operator of the problem
\eqref{eq:lindualproblem}.
\begin{theorem}\label{th:lintruncprobsolvable}
Given $\kappa_0>0$ with $\kappa_0 \ge 1$ if $d=2$,
assume that $\kappa \ge \kappa_0$ and $R \ge R_0$.
In addition, $\kappa_0\ge 1$ is required for $d=2$.
If
\begin{equation}\label{eq:adjapproxpropcond}
\eta(V_h,L_2(B_R))
:=\sup_{f\in L_2(B_R)}\inf_{v_h\in V_h}\frac{\|(\mathcal{A}_N^*)^{-1}f - v_h\|_V}{\|f\|_{0,2,B_R}}
\le \frac{C_-}{2\kappa C_\mathrm{b}},
\end{equation}
where $C_\mathrm{b}>0$ is the continuity constant of the sesquilinear form $a_N$,
then, for any antilinear continuous functional $\ell:\; V\to\C$,
there exists a constant $N^\ast > 0$ such that for $N \ge N^\ast$ the problem

Find $u_h\in V_h$ such that
\begin{equation}\label{eq:lintruncfemproblem}
a_N(u_h,v_h) = \ell(v_h)
\quad\text{for all } v_h\in V_h
\end{equation}
is uniquely solvable.
\end{theorem}
The quantity $\eta(V_h,L_2(B_R))$
is a measure of the approximation ability of $V_h$ to the solutions
of the adjoint problem \eqref{eq:lindualproblem}
over a class of right-hand sides (here $L_2(B_R)$) \cite{Sauter:06}.
\begin{proof}
First we show that the problem \eqref{eq:lintruncfemproblem} has at most one solution.

For this we first \emph{assume} that there exists a solution $u_h\in V_h$ of \eqref{eq:lintruncfemproblem}
and derive an a priori estimate of the error $\|u_N-u_h\|_V$, where $u_N\in V$ is the
solution of the linear problem
\begin{equation}\label{eq:lintruncproblem}
a_N(u_N,v) = \ell(v)
\quad\text{for all } v\in V.
\end{equation}
It is known \cite[Thm.~24]{Angermann:23f} that, under the assumptions of the theorem,
such a solution exists and is unique.

Since $a_N$ satisfies a G{\aa}rding's inequality on $V$ \cite[Lemma~21(ii)]{Angermann:23f},
we have, making use of \eqref{eq:wavenumbernormequiv},
\begin{gather*}
C_-^2 \|u_N-u_h\|_V^2 - 2\kappa^2 \|u_N-u_h\|_{0,2,B_R}^2
\le \Re a_N(u_N-u_h,u_N-u_h)\\
\le |a_N(u_N-u_h,u_N-u_h)|
= |a_N(u_N-u_h,u_N-v_h)|
\quad\text{for all }
v_h\in V_h,
\end{gather*}
where the last relation is a consequence of the Galerkin orthogonality
of the error $u_N-u_h$ of \eqref{eq:lintruncfemproblem}:
\[
a_N(u_N-u_h,v_h) = 0
\quad\text{for all }
v_h\in V_h.
\]
By \cite[Lemma~21(i)]{Angermann:23f}, the sesquilinear form $a_N$ is continuous,
i.e., there exists a constant $C_\mathrm{b} > 0$ independent of $N$ such that
\begin{equation}\label{eq:aNcont}
|a_N(w,v)| \le C_\mathrm{b}\|w\|_V \|v\|_V
\quad\text{for all } w,v\in V.
\end{equation}
Therefore
\begin{equation}\label{eq:intermedGarding}
C_-^2 \|u_N-u_h\|_V^2 - 2\kappa^2 \|u_N-u_h\|_{0,2,B_R}^2
\le C_\mathrm{b}\|u_N-u_h\|_V \|u_N-v_h\|_V
\quad\text{for all }
v_h\in V_h.
\end{equation}
Now we consider the auxiliary adjoint problem \eqref{eq:lindualproblem}
with the particular right-hand side $f:=u_N-u_h$ (keeping the notation $w_N$ for its solution).
Then, again thanks to the Galerkin orthogonality and the continuity of $a_N$,
it follows from \eqref{eq:lindualproblem} that
\begin{align*}
\|u_N-u_h\|_{0,2,B_R}^2
&= \overline{a_N(u_N-u_h,w_N)}
= \overline{a_N(u_N-u_h,w_N-v_h)}\\
&\le C_\mathrm{b}\|u_N-u_h\|_V \|w_N-v_h\|_V
\quad\text{for all }
v_h\in V_h.
\end{align*}
Hence, if $\|u_N-u_h\|_{0,2,B_R}\ne 0$,
\[
\|u_N-u_h\|_{0,2,B_R}\le C_\mathrm{b}\|u_N-u_h\|_V \frac{\|w_N-v_h\|_V}{\|u_N-u_h\|_{0,2,B_R}}.
\]
By means of the solution operator of the problem \eqref{eq:lindualproblem}
we can rewrite this estimate as
\[
\|u_N-u_h\|_{0,2,B_R}
\le C_\mathrm{b}\|u_N-u_h\|_V \frac{\|(\mathcal{A}_N^*)^{-1}(u_N-u_h)-v_h\|_V}{\|u_N-u_h\|_{0,2,B_R}}.
\]
Since the left-hand side does not depend on $v_h\in V_h$, we obtain
\[
\|u_N-u_h\|_{0,2,B_R}
\le C_\mathrm{b}\|u_N-u_h\|_V \inf_{v_h\in V_h}\frac{\|(\mathcal{A}_N^*)^{-1}(u_N-u_h)-v_h\|_V}{\|u_N-u_h\|_{0,2,B_R}},
\]
and from this we conclude, by the definition \eqref{eq:adjapproxpropcond} of $\eta$, that
\[
\|u_N-u_h\|_{0,2,B_R}
\le C_\mathrm{b}\eta(V_h,L_2(B_R))\|u_N-u_h\|_V.
\]
Clearly this inequality is true also for $\|u_N-u_h\|_{0,2,B_R} = 0$
so that we can remove this interim assumption.

Then it follows from \eqref{eq:intermedGarding} that
\[
\left[C_-^2 - 2\kappa^2C_\mathrm{b}^2\eta^2(V_h,L_2(B_R))\right]\|u_N-u_h\|_V
\le C_\mathrm{b} \|u_N-v_h\|_V
\quad\text{for all }
v_h\in V_h.
\]
Consequently, if the condition \eqref{eq:adjapproxpropcond} is satisfied, we get
\begin{equation}\label{eq:schatzestfem}
\|u_N-u_h\|_V
\le 2\frac{C_\mathrm{b}}{C_-^2} \|u_N-v_h\|_V
\quad\text{for all }
v_h\in V_h.
\end{equation}
Now we are prepared to demonstrate the claim of the theorem.
Since the discrete linear problem \eqref{eq:lintruncfemproblem}
is equivalent to a quadratic system of linear algebraic equations,
it is sufficient to show that the corresponding homogeneous problem
has the trivial solution only.
It can be regarded as an approximate problem to the homogeneous
linear problem \eqref{eq:lintruncproblem} (i.e., $\ell=0$).
Since \eqref{eq:lintruncproblem}
is uniquely solvable \cite[Thm.~24]{Angermann:23f}, the corresponding homogeneous
problem has the only solution $u_N=0$.
Now we assume that the homegeneous linear discrete problem \eqref{eq:lintruncfemproblem}
has a nontrivial solution $u_h\in V_h$.
Then the above estimate \eqref{eq:schatzestfem} implies
\[
\|u_h\|_V
\le 2\frac{C_\mathrm{b}}{C_-^2} \|v_h\|_V
\quad\text{for all }
v_h\in V_h,
\]
and this leads to $\|u_h\|_V=0$ in contrast to the assumption.
\end{proof}

Of course the satisfiability of the condition \eqref{eq:adjapproxpropcond}
is crucial. We will postpone this discussion to the end of this section
(Subsection~\ref{subsec:adjapproxpropcond}) and continue examining
the properties of the finite element method.

\begin{corollary}\label{cor:truncfemstab}
Under the assumptions of Thm.~\ref{th:lintruncprobsolvable},
there exists a number $N^\ast\in\N$ and a constant $C>0$
independent of $N\ge N^\ast$ and $h$ such that the solution of
\eqref{eq:lintruncfemproblem} satisfies the estimate
\[
\|u_h\|_V \le C\|\ell\|_{V^\ast}.
\]
\end{corollary}
\begin{proof}
We return to the proof of Thm.~\ref{th:lintruncprobsolvable} and
point out that the estimate \eqref{eq:schatzestfem} is valid for
solutions $u_N\in V$, $u_h\in V_h$ of the general linear problems \eqref{eq:lintruncproblem}
and \eqref{eq:lintruncfemproblem}, respectively.

Hence, by the triangle inequality,
\begin{align*}
\|u_h\|_V &\le \|u_N\|_V + \|u_N-u_h\|_V
\le \|u_N\|_V + 2\frac{C_\mathrm{b}}{C_-^2} \|u_N-v_h\|_V
\quad\text{for all }
v_h\in V_h.
\end{align*}
The particular choice $v_h:=0$ yields
\[
\|u_h\|_V \le \left(1 + 2\frac{C_\mathrm{b}}{C_-^2}\right) \|u_N\|_V.
\]
Since the estimate
\[
\|u_N\|_V \le \frac{4}{C_- \beta(R,\kappa)}\|\ell\|_{V^\ast}
\]
is known \cite[Lemma~25]{Angermann:23f}, we obtain
\[
\|u_h\|_V
\le \left(1 + 2\frac{C_\mathrm{b}}{C_-^2}\right) \frac{4}{C_- \beta(R,\kappa)}\|\ell\|_{V^\ast}
\le \left(1 + 2\frac{C_\mathrm{b}}{C_-^2}\right) \frac{4}{C_- \beta(R,\kappa)}
\sup_{v\in V\setminus\{0\}}\frac{|a_N(u_N,v)|}{\|v\|_{V,\kappa}}.
\]
\end{proof}

Next we demonstrate that the truncated sesquilinear form $a_N$ satisfies
a \emph{discrete} inf-sup condition with a constant that is also uniform
with respect to the mesh parameter $h$,
at least for a sufficiently large truncation parameter $N$.
\begin{lemma}\label{l:truncinfsup}
Under the assumptions of Thm.~\ref{th:lintruncprobsolvable},
there exists a number $N^\ast\in\N$ such that
\[
\beta_{N^\ast,h}(R,\kappa) := \inf_{w_h\in V_h\setminus\{0\}}\sup_{v_h\in V_h\setminus\{0\}}
\frac{|a_N(w_h,v_h)|}{\|w_h\|_{V,\kappa} \|v_h\|_{V,\kappa}}>0
\]
is independent of $N\ge N^\ast$ and $h$.
\end{lemma}
\begin{proof}
We start by observing that, under the assumptions of the lemma,
the truncated sesquilinear form $a_N$ satisfies the inf-sup condition \eqref{eq:infsuptrunc}
with a constant $\beta_{N^\ast}(R,\kappa)$ independent of $N\ge N^\ast$.
Hence, since $V_h\subset V$,
\[
\beta_{N^\ast}(R,\kappa) \le \inf_{w_h\in V_h\setminus\{0\}}\sup_{v\in V\setminus\{0\}}
\frac{|a_N(w_h,v)|}{\|w_h\|_{V,\kappa} \|v\|_{V,\kappa}}.
\]
Unfortunately, this estimate is unsatisfactory because it is not symmetrical
with respect to the spaces.
To overcome this difficulty, we apply an argument known as \emph{Fortin's trick}
\cite{Fortin:77}.

We define an operator $R_h:\;V\to V_h$ as follows:

Given $v\in V$, find $R_hv\in V_h$ such that
\begin{equation}\label{eq:deffortin}
a_N(w_h,R_hv) = a_N(w_h,v)
\quad\text{for all }
w_h\in V_h.
\end{equation}
Since this problem has a unique solution (cf.\ the last part of the proof
of Thm.~\ref{th:lintruncprobsolvable}), the operator $R_h$ is well-defined.
Moreover, the following estimate holds:
\begin{equation}\label{eq:estfortin}
\|R_hv\|_{V,\kappa} \le C_\mathrm{F}\|v\|_{V,\kappa}
\quad\text{for all }
v\in V,
\end{equation}
where the constant $C_\mathrm{F}>0$ does not depend on $N\ge N^\ast$ and $h$.
It can be proved similar to the estimate from Corollary~\ref{cor:truncfemstab},
since the basic properties of the truncated sesquilinear form $a_N$,
namely the continuity estimate and G{\aa}rding's inequality \cite[Lemma~21]{Angermann:23f},
are symmetrical with respect to both arguments.
Then
\begin{align*}
\sup_{v_h\in V_h\setminus\{0\}}
\frac{|a_N(w_h,v_h)|}{\|v_h\|_{V,\kappa}}
&\ge \sup_{v\in V\setminus\{0\}}
\frac{|a_N(w_h,R_hv)|}{\|R_hv\|_{V,\kappa}}\\
&= \sup_{v\in V\setminus\{0\}}
\frac{|a_N(w_h,v)|}{\|R_hv\|_{V,\kappa}}
\ge \frac{1}{C_\mathrm{F}} \sup_{v\in V\setminus\{0\}}
\frac{|a_N(w_h,v)|}{\|v\|_{V,\kappa}}.
\end{align*}
Consequently, $\beta_{N^\ast,h}(R,\kappa)=\beta_{N^\ast}(R,\kappa) C_\mathrm{F}^{-1}$.
\end{proof}

Analogously to the continuous case, a discrete operator
\[
\mathcal{A}_{N,h}:\; V_h \to V_h^\ast
\]
can be defined via
\[
\mathcal{A}_{N,h}w_h(v_h) := a_N(w_h,v_h)
\quad\text{for all }
w_h,v_h \in V_h.
\]
Then, for given $\ell_h\in V_h^\ast$, the problem

Find $u_h\in V_h$ such that
\[
a_N(u_h,v_h) = \ell_h(v_h)
\quad\text{for all } v_h\in V_h
\]
is equivalent to finding a solution $u_h\in V_h$ of the operator equation
\begin{equation}\label{eq:lintruncfemproblemh}
\mathcal{A}_{N,h}u_h = \ell_h.
\end{equation}

Lemma~\ref{l:truncinfsup} ensures that the problem \eqref{eq:lintruncfemproblemh} is well-posed.
Hence there exists a uniformly bounded solution operator $(\mathcal{A}_{N,h})^{-1}:\;V_h^\ast\to V_h$
of the discrete linear problem \eqref{eq:lintruncfemproblem}.
If we consider only elements $\ell_h\in V_h^\ast$ that are restrictions of elements $\ell\in V^\ast$
to $V_h$, Corollary~\ref{cor:truncfemstab} implies that
\[
\|(\mathcal{A}_{N,h})^{-1}\ell\|_{V,\kappa}
\le \frac{4(C_-^2 + 2C_\mathrm{b})C_+}{C_-^3 \beta(R,\kappa)} \|\ell\|_{V^\ast}.
\]
Furthermore, we define a nonlinear operator $\mathcal{F}_N:\; V \to V^\ast$ by
\begin{align*}
\mathcal{F}_N(w)(v)
&:= n_N(w,v)
\qquad\text{for all }
w,v\in V,
\end{align*}
cf.\ \eqref{def:n_N}.
Then, under suitable conditions to $c,f$ and $u^\mathrm{inc}$
(e.g.\ as given in the subsequent theorem),
the problem \eqref{eq:intscalfemproblem4} is equivalent
to the fixed-point problem
\begin{equation}\label{eq:scalfpptruncatedfrm}
u_h = \mathcal{A}_{N,h}^{-1}\mathcal{F}_N(u_h)
\quad\text{in } V_h
\end{equation}
(note that $V_h\subset V$).
\begin{theorem}\label{th:truncscale_eu_fem}
Under the assumptions of Thm.~\ref{th:lintruncprobsolvable},
let the functions $c$ and $f$ generate locally Lipschitz continuous
Nemycki operators in $V$ and assume that there exist functions $w_f,w_c\in V$ such that
$f(\cdot,w_f)\in L_{p_f/(p_f-1)}(\Om)$ and $c(\cdot,w_c)\in L_{p_c/(p_c-2)}(\Om)$,
respectively, where
$p_c,p_f\in\begin{cases}
[2,\infty), & d=2,\\
[2,6], & d=3,
\end{cases}$.
Furthermore let
$u^\mathrm{inc}\in H^1_\mathrm{loc}(\Om^+)$ be such that
additionally $\Delta u^\mathrm{inc}\in L_{2,\mathrm{loc}}(\Om^+)$ holds.

If there exist numbers $\varrho>0$ and $L_\mathcal{F}\in (0,\beta_{N^\ast,h}(R,\kappa))$
(where $N^\ast$ and $\beta_{N^\ast,h}(R,\kappa)$ are from Lemma~\ref{l:truncinfsup})
such that the
the conditions \eqref{eq:thtrunceu1}, \eqref{eq:thtrunceu2}
are satisfied for all $w,v\in K_\varrho^\mathrm{cl}$,
then the problem \eqref{eq:intscalfemproblem4} has a unique solution
$u_h \in K_\varrho^\mathrm{cl}\cap V_h$ for all $N \ge N^\ast$.
\end{theorem}

\begin{remark}
Since $V_h\subset V$, the norm of the operator $R_h$ satisfies $C_\mathrm{F}\ge 1$,
and thus $\beta_{N^\ast,h}\le\beta_{N^\ast}$ (cf.\ the end of the proof
of Lemma~\ref{l:truncinfsup}).
\end{remark}

\begin{proof}
First we mention that $K_\varrho^\mathrm{cl}\cap V_h$ is a closed nonempty subset of $V_h$.

Next we show that
$\mathcal{A}_{N,h}^{-1}\mathcal{F}(K_\varrho^\mathrm{cl}) \subset K_\varrho^\mathrm{cl}\cap V_h$.
To this end we make use of the following estimates, which hold even for all $w,v\in V$:
\begin{align*}
\|vf(\cdot,w)\|_{0,1,\Om}
&\le \|vf(\cdot,w_f)\|_{0,1,\Om} + \|v(f(\cdot,w) - f(\cdot,w_f))\|_{0,1,\Om}\\
&\le \|vf(\cdot,w_f)\|_{0,1,\Om} + \|vL_f(\cdot,w,w_f) |w-w_f|\|_{0,1,\Om}\\
&\le \|v\|_{0,p_f,\Om}\|f(\cdot,w_f)\|_{0,\tilde{q}_f,\Om}
+ \|v\|_{0,p_f,\Om}\|L_f(\cdot,w,w_f)\|_{0,q_f,\Om}\|w-w_f\|_{0,p_f,\Om}\\
&\le \left[\|f(\cdot,w_f)\|_{0,\tilde{q}_f,\Om}
+ \|L_f(\cdot,w,w_f)\|_{0,q_f,\Om}(\|w\|_V + \|w_f\|_V)\right]\|v\|_V,\\
\|zvc(\cdot,w)\|_{0,1,\Om}
&\le \|zvc(\cdot,w_c)\|_{0,1,\Om} + \|zv(c(\cdot,w) - c(\cdot,w_c))\|_{0,1,\Om}\\
&\le \|zvc(\cdot,w_c)\|_{0,1,\Om} + \|zvL_c(\cdot,w,w_c) |w-w_c|\|_{0,1,\Om}\\
&\le \|z\|_{0,p_c,\Om}\|v\|_{0,p_c,\Om}\|c(\cdot,w_c)\|_{0,\tilde{q}_c,\Om}\\
&\quad + \|z\|_{0,p_c,\Om}\|v\|_{0,p_c,\Om}\|L_c(\cdot,w,w_c)\|_{0,q_c,\Om}\|w-w_c\|_{0,p_c,\Om}\\
&\le \left[\|c(\cdot,w_c)\|_{0,\tilde{q}_c,\Om}
+ \|L_c(\cdot,w,w_c)\|_{0,q_c,\Om}(\|w\|_V + \|w_c\|_V)\right]\|z\|_V\|v\|_V
\end{align*}
with
$\frac{1}{p_f}+\frac{1}{\tilde{q}_f}=1$
and
$\frac{2}{p_c}+\frac{1}{\tilde{q}_c}=1$.
With that we obtain
\begin{align*}
\|\mathcal{F}(w)\|_{V^\ast}
&\le \|\ell^\mathrm{contr}(w)\|_{V^\ast} + \|\ell^\mathrm{src}(w)\|_{V^\ast}
+\|\ell^\mathrm{inc}\|_{V^\ast}\\
&\le \kappa^2\left[\|c(\cdot,w_c)-1\|_{0,\tilde{q}_c,\Om}
+ \|L_c(\cdot,w,w_c)\|_{0,q_c,\Om}(\|w\|_V + \|w_c\|_V)\right]\|w\|_V\\
&\quad + \left[\|f(\cdot,w_f)\|_{0,\tilde{q}_f,\Om}
+ \|L_f(\cdot,w,w_f)\|_{0,q_f,\Om}(\|w\|_V + \|w_f\|_V)\right]
+ \|\ell^\mathrm{inc}\|_{V^\ast}\\
&\le \kappa^2\left[\|c(\cdot,w_c)-1\|_{0,\tilde{q}_c,\Om}
+ \|L_c(\cdot,w,w_c)\|_{0,q_c,\Om}(\varrho + \|w_c\|_V)\right]\varrho\\
&\quad + \left[\|f(\cdot,w_f)\|_{0,\tilde{q}_f,\Om}
+ \|L_f(\cdot,w,w_f)\|_{0,q_f,\Om}(\varrho + \|w_f\|_V)\right]\\
&\quad
+ C_\mathrm{tr}\|\vecr\cdot\nabla u^\mathrm{inc} - T_\kappa u^\mathrm{inc}\|_{-1/2,2,S_R}\,.
\end{align*}
Hence the relationship \eqref{eq:thtrunceu1} implies $\|\mathcal{A}_{N,h}^{-1}\mathcal{F}(w)\|_V \le \varrho$.

It remains to show that the mapping $\mathcal{A}_{N,h}^{-1}\mathcal{F}$ is a contraction.

We start with the consideration of the contrast term. From the elementary
decomposition
\[
(c(\cdot,w) - 1) w - (c(\cdot,v) - 1) v
=(c(\cdot,w) - c(\cdot,v)) w + (c(\cdot,v)-1) (w - v)
\]
we see that
\begin{align*}
&\|\ell^\mathrm{contr}(w) - \ell^\mathrm{contr}(v)\|_{V^\ast}\\
&\le \kappa^2\|L_c(\cdot,w,v)\|_{0,q_c,\Om}\|w-v\|_V\|w\|_V\\
&\quad + \kappa^2\left[\|c(\cdot,w_c)-1\|_{0,\tilde{q}_c,\Om}
+ \|L_c(\cdot,w,w_c)\|_{0,q_c,\Om}\|w-w_c\|_V\right]\|w-v\|_V\\
&\le \kappa^2\|L_c(\cdot,w,v)\|_{0,q_c,\Om}\|w-v\|_V\varrho\\
&\quad + \kappa^2\left[\|c(\cdot,w_c)-1\|_{0,\tilde{q}_c,\Om}
+ \|L_c(\cdot,w,w_c)\|_{0,q_c,\Om}(\varrho + \|w_c\|_V)\right]\|w-v\|_V\\
&\le \kappa^2\left[\|L_c(\cdot,w,v)\|_{0,q_c,\Om}\varrho
+ \|c(\cdot,w_c)-1\|_{0,\tilde{q}_c,\Om}
+ \|L_c(\cdot,w,w_c)\|_{0,q_c,\Om}(\varrho + \|w_c\|_V)\right]\|w-v\|_V.
\end{align*}
The estimate of the source term follows immediately from the properties of $f$:
\[
\|\ell^\mathrm{src}(w) - \ell^\mathrm{src}(v)\|_{V^\ast}
\le \|L_f(\cdot,w,v)\|_{0,q_f,\Om}\|w-v\|_V.
\]
From
\[
\|\mathcal{F}(w) - \mathcal{F}(v)\|_{V^\ast}
\le \|\ell^\mathrm{contr}(w) - \ell^\mathrm{contr}(v)\|_{V^\ast}
+ \|\ell^\mathrm{src}(w) - \ell^\mathrm{src}(v)\|_{V^\ast}
\]
and relationship \eqref{eq:thtrunceu2} we thus obtain
\[
\|\mathcal{F}(w) - \mathcal{F}(v)\|_{V^\ast}
\le L_\mathcal{F} \|w-v\|_V.
\]
In summary, Banach's fixed point theorem can be applied (see e.g.\ \cite[Sect.~9.2.1]{Evans:15})
and we conclude that the problem \eqref{eq:intscalfemproblem4}
has a unique solution $u_h \in K_\varrho^\mathrm{cl}\cap V_h$.
\end{proof}

\subsection{FEM-convergence}
In this section we show that the numerical solution is a quasi-optimal
approximation to the solution of the truncated problem.
\begin{theorem}\label{th:truncscale_errest_fem}
Under the assumptions of Thm.~\ref{th:lintruncprobsolvable},
let the functions $c$ and $f$ generate locally Lipschitz continuous
Nemycki operators in $V$ and assume that there exist functions $w_f,w_c\in V$ such that
$f(\cdot,w_f)\in L_{p_f/(p_f-1)}(\Om)$ and $c(\cdot,w_c)\in L_{p_c/(p_c-2)}(\Om)$,
respectively, where
$p_c,p_f\in\begin{cases}
[2,\infty), & d=2,\\
[2,6], & d=3,
\end{cases}$.
Furthermore let
$u^\mathrm{inc}\in H^1_\mathrm{loc}(\Om^+)$ be such that
additionally $\Delta u^\mathrm{inc}\in L_{2,\mathrm{loc}}(\Om^+)$ holds.

If there exist numbers $\varrho>0$ and $L_\mathcal{F}\in (0,\beta_{N^\ast,h}(R,\kappa)\min\{1,C_-^2\})$
(where $N^\ast$ and $\beta_{N^\ast,h}(R,\kappa)$ are from Lemma~\ref{l:truncinfsup})
such that the the conditions \eqref{eq:thtrunceu1}, \eqref{eq:thtrunceu2}
are satisfied for all $w,v\in K_\varrho^\mathrm{cl}$,
then the solution of problem \eqref{eq:intscalfemproblem4} satisfies the estimate
\[
\|u_N-u_h\|_V \le \frac{C_-^2\beta_{N^\ast,h}(R,\kappa) + C_\mathrm{b}}{C_-^2\beta_{N^\ast,h}(R,\kappa) - L_\mathcal{F}}
\inf_{w_h\in V_h}\|u_N-w_h\|_V.
\]
\end{theorem}
\begin{remark}
In the case $d=2$, the parameters $\varrho>0$ and $L_\mathcal{F}$ are the same
as in Thm.~\ref{th:truncscale_eu_fem}. If $d=3$, the condition on $L_\mathcal{F}$
is more restrictive for (the less interesting) wave numbers $\kappa < 1$,
cf.\ the definition of $C_-$ following \eqref{eq:wavenumbernormequiv}.
\end{remark}
\begin{proof}
By Lemma~\ref{l:truncinfsup}, it holds that
\[
\beta_{N^\ast,h}(R,\kappa)\|w_h\|_{V,\kappa}
\le \sup_{v_h\in V_h\setminus\{0\}} \frac{|a_N(w_h,v_h)|}{\|v_h\|_{V,\kappa}}
\quad\text{for all }
w_h\in V_h.
\]
Therefore, replacing $w_h$ by $u_h-w_h$, we also have that
\[
\beta_{N^\ast,h}(R,\kappa)\|u_h-w_h\|_{V,\kappa}
\le \sup_{v_h\in V_h\setminus\{0\}} \frac{|a_N(u_h-w_h,v_h)|}{\|v_h\|_{V,\kappa}}
\quad\text{for all }
w_h\in V_h.
\]
In the numerator on the right-hand side we first write
\[
a_N(u_h-w_h,v_h)
= a_N(u_h-u_N,v_h) + a_N(u_N-w_h,v_h).
\]
If we restrict the problem \eqref{eq:intscalproblem4} to test functions from
$V_h\subset V$ and subtract the result from \eqref{eq:intscalfemproblem4}, we get
\[
a_N(u_h - u_N,v_h) = a_N(u_h,v_h) - a_N(u_N,v_h) = n_N(u_h,v_h) - n_N(u_N,v_h)
\quad\text{for all }
v_h\in V_h.
\]
Therefore
\begin{align}
&\qquad\beta_{N^\ast,h}(R,\kappa)\|u_h-w_h\|_{V,\kappa}\nonumber\\
&\le \sup_{v_h\in V_h\setminus\{0\}} \frac{|n_N(u_h,v_h) - n_N(u_N,v_h)|}{\|v_h\|_{V,\kappa}}
+ \sup_{v_h\in V_h\setminus\{0\}} \frac{|a_N(u_N-w_h,v_h)|}{\|v_h\|_{V,\kappa}}
\label{eq:ees1}\\
&\le \sup_{v_h\in V_h\setminus\{0\}} \frac{|n_N(u_h,v_h) - n_N(u_N,v_h)|}{\|v_h\|_{V,\kappa}}
+ \frac{C_\mathrm{b}}{C_-}\|u_N-w_h\|_V
\quad\text{for all }
w_h\in V_h,
\nonumber
\end{align}
where in the last step we have used \eqref{eq:aNcont} together with \eqref{eq:wavenumbernormequiv}.
Since the second term finally turns into a best approximation error,
it essentially remains to estimate the first term of the upper bound.
Since
\begin{align*}
n_N(u_h,v_h) - n_N(u_N,v_h)
&= \langle\ell^\mathrm{contr}(u_h) - \ell^\mathrm{contr}(u_N),v_h\rangle
+ \langle\ell^\mathrm{src}(u_h) - \ell^\mathrm{src}(u_N),v_h\rangle
\end{align*}
(see \eqref{def:n_N}), we can make us of the estimates given in the second part
of the proof of Thm.~\ref{th:truncscale_eu_fem}:
\begin{align*}
\|\ell^\mathrm{contr}(u_h) - \ell^\mathrm{contr}(u_N)\|_{V_h^\ast}
&\le \kappa^2\big[\|L_c(\cdot,u_h,u_N)\|_{0,q_c,\Om}\varrho
+ \|c(\cdot,w_c)-1\|_{0,\tilde{q}_c,\Om}\\
&\quad + \|L_c(\cdot,u_h,w_c)\|_{0,q_c,\Om}(\varrho + \|w_c\|_V)\big]\|u_h-u_N\|_V.
\end{align*}
The estimate of the source term follows immediately from the properties of $f$:
\[
\|\ell^\mathrm{src}(u_h) - \ell^\mathrm{src}(u_N)\|_{V_h^\ast}
\le \|L_f(\cdot,u_h,u_N)\|_{0,q_f,\Om}\|u_h-u_N\|_V.
\]
From the relationship \eqref{eq:thtrunceu2} we thus obtain
\begin{align*}
\sup_{v_h\in V_h\setminus\{0\}} \frac{|n_N(u_h,v_h) - n_N(u_N,v_h)|}{\|v_h\|_{V,\kappa}}
&\le \frac{1}{C_-}\sup_{v_h\in V_h\setminus\{0\}} \frac{|n_N(u_h,v_h) - n_N(u_N,v_h)|}{\|v_h\|_V}\\
&\le \frac{L_\mathcal{F}}{C_-} \|u_h-u_N\|_V.
\end{align*}
Using this estimate in \eqref{eq:ees1} together with \eqref{eq:wavenumbernormequiv}, we arrive at
\begin{align*}
C_-\beta_{N^\ast,h}(R,\kappa)\|u_h-w_h\|_V
&\le \frac{L_\mathcal{F}}{C_-} \|u_h-u_N\|_V + \frac{C_\mathrm{b}}{C_-}\|u_N-w_h\|_V \\
&\le \frac{L_\mathcal{F}}{C_-} \|u_h-w_h\|_V + \frac{L_\mathcal{F} + C_\mathrm{b}}{C_-}\|u_N-w_h\|_V.
\end{align*}
Hence
\[
\|u_h-w_h\|_V
\le \frac{L_\mathcal{F} + C_\mathrm{b}}{C_-^2\beta_{N^\ast,h}(R,\kappa) - L_\mathcal{F}}\|u_N-w_h\|_V.
\]
Finally, the triangle inequality yields
\[
\|u_N-u_h\|_V \le \|u_N-w_h\|_V + \|u_h-w_h\|_V
\le \left(1 + \frac{L_\mathcal{F} + C_\mathrm{b}}{C_-^2\beta_{N^\ast,h}(R,\kappa) - L_\mathcal{F}}\right)
\|u_N-w_h\|_V.
\]
\end{proof}

\subsection{Discussion of the condition \eqref{eq:adjapproxpropcond}}
\label{subsec:adjapproxpropcond}

First we note that a simple density argument $\bigcup_{h\in\mathcal{H}} V_h \subset V$
is not applicable since we need a uniform estimate w.r.t.\ $f$ in $(\mathcal{A}_N^*)^{-1}f$.
Therefore we have to investigate the adjoint problem \eqref{eq:lindualproblem} more closely.
From
\[
a_N(v,w_N) = (\nabla v,\nabla w_N)_{B_R} - \kappa^2 (v,w_N)_{B_R} - (T_{\kappa,N} v,w_N)_{S_R}
\]
it can be seen immediately that
\begin{align*}
\overline{a_N(v,w_N)} &= (\nabla w_N,\nabla v)_{B_R} - \kappa^2 (w_N,v)_{B_R} - (w_N,T_{\kappa,N} v)_{S_R}\\
&= (\nabla w_N,\nabla v)_{B_R} - \kappa^2 (w_N,v)_{B_R} - (T_{\kappa,N}^\ast w_N, v)_{S_R},
\end{align*}
where $T_{\kappa,N}^\ast$ is the adjoint operator of $T_{\kappa,N}$ in $L_2(S_R)$.
So the structure of the adjoint problem \eqref{eq:lindualproblem} differs from that
of problem \eqref{eq:lintruncproblem} only by a different DtN term.

Hence the question of a representation of the adjoint operator and its properties arises.
The answer is based on the observation that, if $w:\;B_R^+ \to\C$ is a solution
of the exterior Dirichlet problem
\begin{equation}\label{eq:extDirproblem}
\begin{aligned}
-\Delta w - \kappa^2 w = 0
\quad\text{in }
B_R^+,\\
w = g
\quad\text{on }
S_R,\\
\lim_{|\vecx|\to\infty} |\vecx|^{(d-1)/2}
\left(\vecr\cdot\nabla w - i\kappa w\right) = 0,
\end{aligned}
\end{equation}
where $g:\;S_R\to\C$ is given, then the conjugate $\overline{w}$ is a solution
of an analogous exterior Dirichlet problem
but with an \emph{ingoing} radiation condition:
\begin{equation}\label{eq:extDirproblemingoing}
\begin{aligned}
-\Delta \overline{w} - \kappa^2 \overline{w} = 0
\quad\text{in }
B_R^+,\\
\overline{w} = \overline{g}
\quad\text{on }
S_R,\\
\lim_{|\vecx|\to\infty} |\vecx|^{(d-1)/2}
\left(\vecr\cdot\nabla \overline{w} + i\kappa \overline{w}\right) = 0.
\end{aligned}
\end{equation}
Therefore, if $T_\kappa$ is the (exact) Dirichlet-to-Neumann operator associated with
\eqref{eq:extDirproblem}, i.e.
\begin{equation}\label{def:smoothDtN}
g\mapsto T_\kappa g := \left.\vecr\cdot\nabla w\right|_{S_R},
\end{equation}
and if $T_\kappa^+$ is the Dirichlet-to-Neumann operator associated with
\eqref{eq:extDirproblemingoing}, then
\begin{equation}\label{eq:DtNrelation}
T_\kappa^+ \overline{g} := \left.\vecr\cdot\nabla \overline{w}\right|_{S_R}
= \overline{\left.\vecr\cdot\nabla w\right|_{S_R}}
= \overline{T_\kappa g}.
\end{equation}
Therefore
\[
T_\kappa^+ v = \overline{T_\kappa \overline{v}}.
\]
According to \cite[p.~72]{Colton:13b} it is possible to develop a theory of incident waves
that fulfill the last condition in \eqref{eq:extDirproblemingoing} analogous to the previous
case of outgoing waves that satisfy Sommerfeld's radiation condition in \eqref{eq:extDirproblem}.
Therefore the corresponding DtN operator also has the analogous properties,
especially with regard to regularity, cf., eg.~\cite[Thm.~2]{Angermann:23f}.

Since only the truncated operator is required in the adjoint problem \eqref{eq:lindualproblem},
we investigate a representation motivated by \eqref{eq:DtNrelation}.

We start with the two-dimensional situation.
So let
\begin{equation}\label{eq:wvseries2}
\begin{aligned}
v(\vecx) &= v(R\hat\vecx)
= \sum_{|k|\in\N_0} v_k(R)Y_k(\hat\vecx),
\quad
\vecx=R\hat\vecx\in S_R,
\end{aligned}
\end{equation}
be the series representation of $v|_{S_R}$ with the Fourier coefficients
(cf.~\eqref{eq:Fouriercoeff2d})
\begin{align*}
v_k(R) &= (v(R\cdot),Y_k)_{S_1}
= \int_{S_1}v(R\hat\vecx)\overline{Y_k}(\hat\vecx)ds(\hat\vecx)
\end{align*}
(analogously for $w|_{S_R}$).
Then the Fourier representation of $\overline{v}$ formally reads as
\begin{align*}
\overline{v}(\vecx) &= \overline{v}(R\hat\vecx)
= \sum_{|k|\in\N_0} \tilde{v}_k(R)Y_k(\hat\vecx),
\quad
\vecx=R\hat\vecx\in S_R,
\end{align*}
where
\begin{align*}
\tilde{v}_k(R) &= (\overline{v}(R\cdot),Y_k)_{S_1}
= \int_{S_1}\overline{v}(R\hat\vecx)\overline{Y_k}(\hat\vecx)ds(\hat\vecx).
\end{align*}
On the other hand, using \eqref{eq:wvseries2},
\begin{align*}
\overline{v(\vecx)} &= \overline{v(R\hat\vecx)}
= \sum_{|k|\in\N_0} \overline{v_k(R)}\,\overline{Y_k(\hat\vecx)}
= \sum_{|k|\in\N_0} \overline{v_k(R)}Y_{-k}(\hat\vecx)
= \sum_{|k|\in\N_0} \overline{v_{-k}(R)}Y_k(\hat\vecx),
\quad
\vecx=R\hat\vecx\in S_R.
\end{align*}
Hence, by uniqueness of the Fourier representation,
$\tilde{v}_k(R) = \overline{v_{-k}(R)}$.
Then, by \eqref{def:2dDtNtrunc},
\begin{align*}
T_{\kappa,N} \overline{v}(\vecx)
&= \frac{1}{R}\sum_{|k|\le N}
Z_k(\kappa R)\tilde{v}_k(R)Y_k(\hat\vecx)\\
&= \frac{1}{R}\sum_{|k|\le N}
Z_k(\kappa R)\overline{v_{-k}(R)}Y_k(\hat\vecx),
\quad
\vecx=R\hat\vecx\in S_R,
\end{align*}
and
\begin{align*}
\overline{T_{\kappa,N} \overline{v}(\vecx)}
&= \frac{1}{R}\sum_{|k|\le N}
\overline{Z_k(\kappa R)}v_{-k}(R)\overline{Y_k(\hat\vecx)}
= \frac{1}{R}\sum_{|k|\le N}
\overline{Z_k(\kappa R)}v_{-k}(R)Y_{-k}(\hat\vecx)\\
&= \frac{1}{R}\sum_{|k|\le N}
\overline{Z_{-k}(\kappa R)}v_k(R)Y_k(\hat\vecx),
\quad
\vecx=R\hat\vecx\in S_R.
\end{align*}
Thus, by an analogous relation to \eqref{eq:DtNrelation} and
the orthonormality of the circular harmonics \cite[Prop.~3.2.1]{Zeidler:95as},
\begin{align*}
\left(w,T_{\kappa,N}^+ v\right)_{S_R}
&= \left(w,\overline{T_{\kappa,N} \overline{v}}\right)_{S_R}
= \frac{1}{R}\sum_{|n|,|k|\le N}
\left(w_n(R)Y_n(R^{-1}\cdot),\overline{Z_{-k}(\kappa R)}v_k(R)Y_k(R^{-1}\cdot)\right)_{S_R}\\
&= \sum_{|n|,|k|\le N}
Z_{-k}(\kappa R)\left(w_n(R)Y_n,v_k(R)Y_k\right)_{S_1}\\
&= \sum_{|n|\le N}
Z_{-n}(\kappa R)\left(w_n(R)Y_n,v_n(R)Y_n\right)_{S_1}.
\end{align*}
Since $H^{(1)}_{-n}(\xi)=(-1)^nH^{(1)}_{n}(\xi)$
\cite[(10.4.2)]{NIST:23}, we have $Z_{-n}=Z_n$ and get
\begin{align*}
\left(w,T_{\kappa,N}^+ v\right)_{S_R}
&= \sum_{|n|\le N} Z_n(\kappa R)\left(w_n(R)Y_n,v_n(R)Y_n\right)_{S_1}\\
&= \sum_{|n|,|k|\le N}
Z_n(\kappa R)\left(w_n(R)Y_n,v_k(R)Y_k\right)_{S_1}\\
&= \frac{1}{R}\sum_{|n|,|k|\le N}
\left(Z_n(\kappa R) w_n(R)Y_n(R^{-1}\cdot),v_k(R)Y_k(R^{-1}\cdot)\right)_{S_R}
= \left(T_{\kappa,N}w,v\right)_{S_R}.
\end{align*}
This shows that the operator $T_{\kappa,N}^+$ is the adjoint of $T_{\kappa,N}$.

\medskip
The investigation of the case $d=3$ runs similarly. So let
\begin{equation}\label{eq:wvseries3}
\begin{aligned}
v(\vecx) &= v(R\hat\vecx)
= \sum_{k\in\N_0}\sum_{|l|\le k} v_k^l(R)Y_k^l(\hat\vecx),
\quad
\vecx=R\hat\vecx\in S_R,
\end{aligned}
\end{equation}
be the series representation of $v|_{S_R}$ with the Fourier coefficients
(cf.~\eqref{eq:Fouriercoeff3d})
\begin{align*}
v_k^l(R) &= (v(R\cdot),Y_k^l)_{S_1}
= \int_{S_1}v(R\hat\vecx)\overline{Y_k^l}(\hat\vecx)ds(\hat\vecx)
\end{align*}
(analogously for $w|_{S_R}$).
Then the Fourier representation of $\overline{v}$ formally reads as
\begin{align*}
\overline{v}(\vecx) &= \overline{v}(R\hat\vecx)
= \sum_{k\in\N_0}\sum_{|l|\le k} \tilde{v}_k^l(R)Y_k^l(\hat\vecx),
\quad
\vecx=R\hat\vecx\in S_R,
\end{align*}
where
\begin{align*}
\tilde{v}_k^l(R) &= (\overline{v}(R\cdot),Y_k^l)_{S_1}
= \int_{S_1}\overline{v}(R\hat\vecx)\overline{Y_k^l}(\hat\vecx)ds(\hat\vecx).
\end{align*}
On the other hand, using \eqref{eq:wvseries3} and the fact that the
associated Legendre functions of the first kind $P_k^l$
are real-valued for real-valued arguments \cite[Sect.~14.3, 15.2]{NIST:23},
we have that
\begin{align*}
\overline{v(\vecx)} &= \overline{v(R\hat\vecx)}
= \sum_{k\in\N_0}\sum_{|l|\le k} \overline{v_k^l(R)}\,\overline{Y_k^l(\hat\vecx)}
= \sum_{k\in\N_0}\sum_{|l|\le k} \overline{v_k^l(R)}Y_k^{-l}(\hat\vecx)\\
&= \sum_{k\in\N_0}\sum_{|l|\le k} \overline{v_k^{-l}(R)}Y_k^l(\hat\vecx),
\quad
\vecx=R\hat\vecx\in S_R.
\end{align*}
Hence, by uniqueness of the Fourier representation,
$\ \tilde{v}_k^l(R) = \overline{v_k^{-l}(R)}$.
Then, by \eqref{def:3dDtNtrunc},
\begin{align*}
T_{\kappa,N} \overline{v}(\vecx)
&:= \frac{1}{R}\sum_{k=0}^N\sum_{|l|\le k}
z_k(\kappa R)\tilde{v}_k^l(R)Y_k^l(\hat\vecx)\\
&= \frac{1}{R}\sum_{k=0}^N\sum_{|l|\le k}
z_k(\kappa R)\overline{v_k^{-l}(R)}Y_k^l(\hat\vecx),
\quad
\vecx=R\hat\vecx\in S_R,
\end{align*}
and
\begin{align*}
\overline{T_{\kappa,N} \overline{v}(\vecx)}
&= \frac{1}{R}\sum_{k=0}^N\sum_{|l|\le k}
\overline{z_k(\kappa R)}v_k^{-l}(R)\overline{Y_k^l(\hat\vecx)}\\
&= \frac{1}{R}\sum_{k=0}^N\sum_{|l|\le k}
\overline{z_k(\kappa R)}v_k^{-l}(R)Y_k^{-l}(\hat\vecx)\\
&= \frac{1}{R}\sum_{k=0}^N\sum_{|l|\le k}
\overline{z_k(\kappa R)}v_k^l(R)Y_k^l(\hat\vecx),
\quad
\vecx=R\hat\vecx\in S_R.
\end{align*}
Thus, by an analogous relation to \eqref{eq:DtNrelation} and
the orthonormality of the spherical harmonics \cite[Thm.~2.8]{Colton:19},
\begin{align*}
\left(w,T_{\kappa,N}^+ v\right)_{S_R}
&= \left(w,\overline{T_{\kappa,N} \overline{v}}\right)_{S_R}\\
&= \frac{1}{R}\sum_{n,k\in\N_0}\sum_{|m|\le n,|l|\le k}
\left(w_n^m(R)Y_n^m(R^{-1}\cdot),\overline{z_k(\kappa R)}v_k^l(R)Y_k^l(R^{-1}\cdot)\right)_{S_R}\\
&= R\sum_{n,k\in\N_0}\sum_{|m|\le n,|l|\le k}
z_k(\kappa R)\left(w_n^m(R)Y_n^m,v_k^l(R)Y_k^l\right)_{S_1}\\
&= R\sum_{n\in\N_0}\sum_{|m|\le n}
z_n(\kappa R) w_n^m(R)\overline{v_n^m}(R)\\
&= R\sum_{n,k\in\N_0}\sum_{|m|\le n,|l|\le k}
z_n(\kappa R)\left(w_n^m(R)Y_n^m,v_k^l(R)Y_k^l\right)_{S_1}\\
&= \frac{1}{R}\sum_{n,k\in\N_0}\sum_{|m|\le n,|l|\le k}
\left(z_n(\kappa R)w_n^m(R)Y_n^m(R^{-1}\cdot),v_k^l(R)Y_k^l(R^{-1}\cdot)\right)_{S_R}
= \left(T_{\kappa,N}w,v\right)_{S_R}.
\end{align*}
This shows that the operator $T_{\kappa,N}^+$ is the adjoint of $T_{\kappa,N}$.

In preparation for the concluding discussion of the quantity $\eta(V_h,L_2(B_R))$
we prove the following lemma, which is a generalization of \cite[Lemma~22]{Angermann:23f}.

\begin{lemma}\label{l:bilintruncerr}
For given $w\in H^s(S_R)$, $s\ge 1/2$, and $v\in H^{1/2}(S_R)$ it holds that
\[
\left|\left((T_\kappa - T_{\kappa,N}) w,v\right)_{S_R}\right| \le c(N,R,s)\|w\|_{s,2,S_R}\|v\|_{1/2,2,S_R},
\]
where \quad $\displaystyle c(N,R,s) := \frac{d-1 + |\kappa R|^2)^{1/2}}{R(1+N^2)^{(2s-1)/4}}$\,.
\end{lemma}
\begin{proof}
We start with the two-dimensional situation and the series representations
\eqref{eq:wvseries2} of $w|_{S_R},v|_{S_R}$.
The norm on the Sobolev space $H^s(S_R)$, $s\ge0$,
can be defined as follows \cite[Ch.~1, Rem.~7.6]{Lions:72a}:
\begin{equation}\label{def:SobFourierNorm2}
\|v\|_{s,2,S_R}^2:=
R\sum_{n\in\Z} (1+n^2)^s |v_n(R)|^2.
\end{equation}
Then, by \eqref{def:2dDtNtrunc},
the orthonormality of the circular harmonics \cite[Prop.~3.2.1]{Zeidler:95as}
and \eqref{def:SobFourierNorm2},
\begin{align*}
\left|\left((T_\kappa - T_{\kappa,N}) w,v\right)_{S_R}\right|
&= \left|\sum_{|n|> N}
Z_n(\kappa R)w_n(R)\overline{v_n}(R)\right|\\
&= \left|\sum_{|n|> N}
\frac{Z_n(\kappa R)}{(1+n^2)^{(2s+1)/4}}(1+n^2)^{s/2}w_n(R)(1+n^2)^{1/4}\overline{v_n}(R)\right|\\
&\le \max_{|n|> N}
\left|\frac{Z_n(\kappa R)}{(1+n^2)^{(2s+1)/4}}\right|
\sum_{|n|> N}\left|(1+n^2)^{s/2}w_n(R)(1+n^2)^{1/4}\overline{v_n}(R)\right|\\
&\le \max_{|n|> N}
\left|\frac{Z_n(\kappa R)}{(1+n^2)^{(2s+1)/4}}\right|
\left(\sum_{|n|> N}(1+n^2)^{1/2}\left|w_n(R)\right|^2\right)^{1/2}\\
&\quad \times
\left(\sum_{|n|> N}(1+n^2)^{1/2}\left|v_n(R)\right|^2\right)^{1/2}\\
&\le \frac{1}{R}\max_{|n|> N}
\left|\frac{Z_n(\kappa R)}{(1+n^2)^{(2s+1)/4}}\right|\|w\|_{s,2,S_R}\|v\|_{1/2,2,S_R}.
\end{align*}
The estimate
\[
\frac{1}{1+n^2} |Z_n(\kappa R)|^2 \le 1 + |\kappa R|^2,
\quad
|n|\in\N,
\]
(see, e.g., \cite[Cor.~5]{Angermann:23f}) implies
\[
\frac{|Z_n(\kappa R)|}{(1+n^2)^{(2s+1)/4}} \le \frac{1 + |\kappa R|^2)^{1/2}}{(1+N^2)^{(2s-1)/4}}
\quad\text{for all }
|n|>N.
\]
The investigation of the case $d=3$ runs similarly, i.e., we start
from the series representations \eqref{eq:wvseries3} of $w|_{S_R},v|_{S_R}$.
The norm on the Sobolev space $H^s(S_R)$, $s\ge0$,
can be defined as follows \cite[Ch.~1, Rem.~7.6]{Lions:72a}:
\begin{equation}\label{def:SobFourierNorm3}
\|v\|_{s,2,S_R}^2:=
R^2\sum_{n\in\N_0}\sum_{|m|\le n} (1+n^2)^s |v_n^m(R)|^2.
\end{equation}
Then, by \eqref{def:3dDtNtrunc},
the orthonormality of the spherical harmonics \cite[Thm.~2.8]{Colton:19}
and \eqref{def:SobFourierNorm3},
\begin{align*}
\left|\left((T_\kappa - T_{\kappa,N}) w,v\right)_{S_R}\right|
&= R\left|\sum_{n> N}\sum_{|m|\le n}
z_n(\kappa R) w_n^m(R)\overline{v_n^m}(R)\right|\\
&= R\left|\sum_{n> N}\sum_{|m|\le n}
\frac{z_n(\kappa R)}{(1+n^2)^{(2s+1)/4}}(1+n^2)^{s/2}w_n^m(R)(1+n^2)^{1/4}\overline{v_n^m}(R)\right|\\
&\le R\max_{n> N}
\left|\frac{z_n(\kappa R)}{(1+n^2)^{(2s+1)/4}}\right|
\sum_{n> N}\sum_{|m|\le n}\left|(1+n^2)^{s/2}w_n^m(R)(1+n^2)^{1/4}\overline{v_n^m}(R)\right|\\
&\le R\max_{n> N}
\left|\frac{z_n(\kappa R)}{(1+n^2)^{(2s+1)/4}}\right|
\left(\sum_{n> N}\sum_{|m|\le n}(1+n^2)^{1/2}\left|w_n^m(R)\right|^2\right)^{1/2}\\
&\quad \times
\left(\sum_{n> N}\sum_{|m|\le n}(1+n^2)^{1/2}\left|v_n^m(R)\right|^2\right)^{1/2}\\
&\le \frac{1}{R}\max_{n> N}
\left|\frac{z_n(\kappa R)}{(1+n^2)^{(2s+1)/4}}\right|\|w\|_{s,2,S_R}\|v\|_{1/2,2,S_R}.
\end{align*}
Thanks to the estimate
\[
\frac{1}{1+n^2} |z_n(\kappa R)|^2 \le 2 + |\kappa R|^2,
\quad
n\in\N_0,
\]
(see, e.g., \cite[Cor.~5]{Angermann:23f}) we can set
\[
c(N,R,s) := \frac{2 + |\kappa R|^2)^{1/2}}{R(1+N^2)^{(2s-1)/4}}\,.
\]
\end{proof}

The above considerations give rise to the conclusion that those analytical properties
of the adjoint exact and truncated DtN operators $T_\kappa^+$, $T_{\kappa,N}^+$
that are relevant here are the same as those of the original operators
$T_\kappa$, $T_{\kappa,N}$ (cf.\ \cite[Thm.~2, Lemmata~22, 23]{Angermann:23f}).
In particular, we may expect the solution $w\in V$ of the \emph{untruncated} adjoint problem
\[
\overline{a(v,w)} = (v,f)_{B_R}
\quad\text{for all } v\in V,
\]
where $a$ is the sesquilinear form defined in \eqref{def:n_N}
with $T_{\kappa,N}$ replaced by $T_\kappa$,
to be $H^2(B_R)$-regular in the sense that there exists a constant
$C(R,\kappa_0)>0$ depending only on $R$ and $\kappa_0$ such that the following
estimate holds:
\begin{equation}\label{eq:adjointreg}
\|w\|_{2,2,B_R} \le \frac{C(R,\kappa_0)}{\kappa}\|f\|_{0,2,B_R}
\end{equation}
(cf.\ \cite[Lemma~3.5]{Melenk:10}, where even more subtle estimates of the solution
of the original problem are given).

This makes it possible to transfer the arguments used in the proof
of \cite[Thm.~24]{Angermann:23f} regarding the solution of problem \eqref{eq:lintruncproblem}
with slight modifications to the adjoint problem \eqref{eq:lindualproblem}.

In particular, we have the estimate
\[
C_-^2 \|w-w_N\|_V \le \eta_1 + 2\kappa^2\eta_2 C_-^{-1}C(R,\kappa) C_\mathrm{emb},
\]
where the positive quantity $\eta_1$ can be estimated by the help of Lemma~\ref{l:bilintruncerr} with $s=3/2$
and the trace theorem \cite[Thm.~3.37]{McLean:00}:
\[
\eta_1 \le c\Big(N,R,\frac{3}{2}\Big)C_\mathrm{tr}^2\|w\|_{2,2,B_R}
= \frac{d-1 + |\kappa R|^2)^{1/2}}{R(1+N^2)^{1/2}}\,C_\mathrm{tr}^2\|w\|_{2,2,B_R}.
\]
Regarding the positive quantity $\eta_2$ we can immediately use the estimate
\[
\eta_2 \le \frac{CC_\mathrm{tr}^2\kappa}{(1+N^2)^{1/2}} \|w_N\|_V.
\]
(analogously to the proof of \cite[Thm.~24]{Angermann:23f}).
Summarizing the above estimates, we obtain
\begin{equation}\label{eq:adjointschatzest1}
\|w-w_N\|_V \le \frac{C}{(1+N^2)^{1/2}} \big[\|w\|_{2,2,B_R} + \|w_N\|_V\big]
\end{equation}
with a constant $C>0$ independent of $N$ (but depending on $R$, $\kappa$).

Now we can turn to the estimation of the quantity $\eta(V_h,L_2(B_R))$.
Based on the relation
\begin{align*}
\inf_{v_h\in V_h}\|(\mathcal{A}_N^*)^{-1}f - v_h\|_V
&= \inf_{v_h\in V_h}\|w_N - v_h\|_V
\le \|w_N - I_h w\|_V\\
&\le \|w_N - w\|_V + \|w - I_h w\|_V,
\end{align*}
we can estimate the first term on the very right-hand side by \eqref{eq:adjointschatzest1}
and the second one by means of (D\ref{item:D4}) (or Thm.~\ref{th:interrest}) with $r=1$.
This results in
\[
\inf_{v_h\in V_h}\|(\mathcal{A}_N^*)^{-1}f - v_h\|_V
\le \frac{C}{(1+N^2)^{1/2}} \big[\|w\|_{2,2,B_R} + \|w_N\|_V\big]
+ C_\mathrm{ip} h |w|_{2,2,B_R}.
\]
Furthermore, since under the assumptions of Thm.~\ref{th:lintruncprobsolvable}
w.r.t.\ $\kappa$ and $R$ the sesquilinear form $a_N$ satisfies an inf-sup condition
\cite[Lemma~25]{Angermann:23f}
(where the inf-sup constant of $a_N$ for all $N$ greater than a certain number $N^\ast\in\N$
no longer depends on $N$),
we have the estimate
\[
\|w_N\|_V \le C_\mathrm{N^\ast,adj}(R,\kappa) \|f\|_{0,2,B_R}
\]
with a positive coefficient $C_\mathrm{N^\ast,adj}(R,\kappa)$,
where the latter is independent of $N$ for all $N\ge N^\ast$.
This estimate together with \eqref{eq:adjointreg} leads to
\[
\inf_{v_h\in V_h}\|(\mathcal{A}_N^*)^{-1}f - v_h\|_V
\le C\Big[\frac{1}{(1+N^2)^{1/2}} + h\Big]\|f\|_{0,2,B_R}
\]
with a constant $C>0$ independent of $N$ for all $N\ge N^\ast$
This relation shows that the condition \eqref{eq:adjapproxpropcond} can be satisfied
provided $N^\ast$ is sufficiently large and $h$ is sufficiently small.

\section{Conclusion}
The work described an exemplary finite element approach for approximating
the solution of a nonlinear Helmholtz equation,
which is reduced to a spherical domain using DtN truncation.
For this, solvability, stability and convergence could be shown.
The satisfiability of the approximation property of the finite element space
required for solvability for the solution of a class of adjoint linear problems was discussed.
Together with the estimate of the truncation error from Thm.~\ref{th:truncscalerrest}
an estimate of the total error, i.e., the error between the FE approximation
and the solution of the full-space transmission problem, is available now.
In perspective, the work is to be continued with the transfer to more robust FE methods
in connection with a (not consistently conducted here) discussion
of the dependence of the constants and parameters on the wave number and the radius
of the computational domain, and with an improvement of the numerical treatment
of the interface $\pOm$.

\newcommand{\etalchar}[1]{$^{#1}$}

\end{document}